\documentclass[reqno,final]{amsart}

\usepackage{caption}
\usepackage{amsfonts,latexsym,graphicx,amsmath,amssymb,bbm,enumitem}

\makeindex

\usepackage[color,notref,notcite]{showkeys}

\definecolor{labelkey}{rgb}{0.6,0,1}

\allowdisplaybreaks

\definecolor{labelkey}{rgb}{0.6,0,1}

\usepackage[normalem]{ulem}
\newcounter{corr}
\definecolor{violet}{rgb}{0.580,0.,0.827}
\newcommand{\corr}[3]{\typeout{Warning : a correction remains in page
\thepage}
				\stepcounter{corr}        
				{\color{blue}\ifmmode\text{\,\sout{\ensuremath{#1}}\,}\else\sout{#1}\fi}
        {\color{red}#2}
        {\color{violet} \fbox{\thecorr}#3}}

\newcounter{cst}
\catcode`\@=11
\def\ctel#1{C_{\refstepcounter{cst}\@bsphack
\protected@write\@auxout{}%
           {\string\newlabel{#1}{{\thecst}{\thepage}}}\thecst}}
\catcode`\@=12

\newcounter{cexp}
\catcode`\@=11
\def\terml#1{T_{\refstepcounter{cexp}\@bsphack
\protected@write\@auxout{}%
           {\string\newlabel{#1}{{\thecexp}{\thepage}}}\thecexp}}
\catcode`\@=12

\newcommand{\mathbi}[1]{{\boldsymbol #1}}

\newcommand{\eop}{{\unskip\nobreak\hfil\penalty50
           \hskip2em\hbox{}\nobreak\hfil\mbox{\rule{1ex}{1ex} \qquad}
   \parfillskip=0pt
   \finalhyphendemerits=0\par\medskip}}

\newtheorem{theorem}{Theorem}[section]
\newtheorem{remark}[theorem]{Remark}

\newtheorem{definition}[theorem]{Definition}

\definecolor{shadecolor}{gray}{0.92}
\definecolor{TFFrameColor}{gray}{0.92}
\definecolor{TFTitleColor}{rgb}{0,0,0}

\newcommand{\ba}{\begin{array}{llll}   }
\newcommand{\bac}{\begin{array}{c}}
\newcommand{\bari}{\begin{array}{r}}
\newcommand{\ea}{\end{array}}

\newcommand{\ban}{\begin{array}{llll}}
\newcommand{\ean}{\end{array}}

\newcommand{\be}{\begin{equation}}
\newcommand{\ee}{\end{equation}}

\newcommand{\beqsys }{\beqtab \left \{ \begin{array}{l}}
\newcommand{\eeqsys }{\end{array} \right . \eeqtab }

\newcommand{\benum}{\begin{enumerate}}
\newcommand{\eenum}{\end{enumerate}}

\newcommand{\beqtab}{\begin{eqnarray}} 
\newcommand{\eeqtab}{\end{eqnarray}}

\newcommand{\combineln}[2]{{\rm(\textbf{#1#2})}}

\newcommand{\bfn}{\mathbf{n}}

\newcommand{\conv}{\rightarrow} 

\newcommand{\darcyU}{\mathbf{u}}

\newcommand{\disc}{{\mathcal D}}
\newcommand{\discC}{{\mathcal C}}
\newcommand{\discP}{{\mathcal P}}

\newcommand{\discDarcyU}[1][(n+1)]{\darcyU_{\discP}^{#1}}

\newcommand{\dt}{{\delta\!t}}
\newcommand{\dtDisc}{\dt^{(n+\frac{1}{2})}}

\renewcommand{\div}{{\mathop{\rm div}}}
\newcommand{\divg}{\div}
\newcommand{\proj}{E}

\newcommand{\ICinterp}{\mathcal{I}}

\newcommand{\mesh}{{\mathcal M}}

\newcommand{\N}{\mathbb N}

\renewcommand{\O}{\Omega}

\newcommand{\R}{\mathbb R}

\newcommand{\x}{\mathbi{x}}

\def\PiD{\Pi_\disc}
\def\PiDc{\Pi_{\discC}}
\def\PiDp{\Pi_{\discP}}

\def\gradDc{\nabla_{\discC}}
\def\gradDp{\nabla_{\discP}}

\usepackage{xparse}
\DeclareDocumentCommand{\RPiD}{ O{\disc} O{,0} }{\Pi_{#1}(X_{#1#2})}

\def\Fdof#1{{\bm{\mathcal{F}}(#1,\R)}}
\makeatletter
\def\Fdof{\@ifnextchar[{\@with}{\@without}}
\def\@with[#1]#2{{\bm{\mathcal{F}}(#2;#1)}}
\def\@without#1{{\bm{\mathcal{F}}(#1,\R)}}
\makeatother

\newcommand{\U}{\mathbf{U}}
\newcommand{\Daru}{\mathbf{u}}

\def\K{\mathbf{K}}

\def\RT0{\mathbb{RT}_0}

\begin{document}
	
	\title[]{A combined GDM--ELLAM--MMOC scheme for advection dominated PDEs}
	
\author{Hanz Martin Cheng}
\address{Department of Mathematics and Computer Science, Eindhoven University of Technology, P.O. Box 513, 5600 MB Eindhoven, the Netherlands.
\texttt{h.m.cheng@tue.nl}}
\author{J\'er\^ome Droniou}
\address{School of Mathematics, Monash University, Clayton, Victoria 3800, Australia.
\texttt{jerome.droniou@monash.edu}}
\author{Kim-Ngan Le}
\address{School of Mathematics, Monash University, Clayton, Victoria 3800, Australia.
\texttt{ngan.le@monash.edu}}

	\date{\today}
	
	%
	%
	
	\maketitle
\begin{abstract}
	We propose a combination of the Eulerian Lagrangian Localised Adjoint Method (ELLAM)  and the Modified Method of Characteristics (MMOC) for 
	time-dependent advection-domina\-ted PDEs. The combined scheme, so-called GEM scheme, takes advantages of both ELLAM scheme (mass conservation) and MMOC scheme (easier computations), while at the same time avoids their disadvantages (respectively, harder tracking around the injection regions, and loss of mass). 
	We present a precise analysis of mass conservation properties for these three schemes, and after achieving global mass balance, an adjustment yielding local volume conservation is then proposed. Numerical results for all three schemes are then compared, illustrating the advantages of the GEM scheme.
	A convergence result of the MMOC scheme, motivated by our previous work~\cite{CDL17-convergence-ELLAM}, is provided which can be extended to obtain the convergence of GEM scheme.
\end{abstract}
\section{Models and assumptions} \label{sec:PureAdvectionModel}

\subsection{Introduction}

In this paper, we consider a time-dependent advection-domina\-ted PDE \eqref{eq:advection}, and study some numerical schemes for this equation that are based on characteristic methods. These types of PDEs are encountered in many important fields, such as mathematical models in porous medium flow (e.g. reservoir simulation), and fluid dynamics (e.g. Navier-–Stokes equations). A short summary, which includes most of the commonly used numerical schemes for \eqref{eq:advection}, together with their advantages and disadvantages, have been presented in \cite{EW01-summary-advection-dominated}. 

In particular, our work focuses on two types of numerical schemes based on characteristic methods which are popularly used, namely the Eulerian Lagrangian Localised Adjoint Method (ELLAM)  and the Modified Method of Characteristics (MMOC). The advantages of these schemes lie on the fact that they are based on characteristic methods, and thus capture the advective component of the PDE better than upwinding schemes. Moreover, these schemes are not limited by CFL constraints, and hence large time steps can be taken for numerical simulations. These are usually combined with finite difference (FD), finite element (FE) or finite volume (FV) discretisations, in order to provide a complete numerical scheme for \eqref{eq:advection}. To cite a few examples, the FE--MMOC \cite{DR82-MMOC-main}, FE--ELLAM \cite{CRHE90-ELLAM-main}, and FV--ELLAM \cite{HR98-FV-ELLAM}, have been used to discretise \eqref{eq:advection}. Other variants of the ELLAM, as well as a summary of its properties, have also been presented in \cite{RC-02-overview}. More recent variants of the ELLAM involve the volume corrected characteristics mixed method (VCCMM) \cite{AH06,AW11-stability-monotonicity-implementation}. Aside from the global mass conservation property of ELLAM, these ensure that local volume conservation is achieved. On the other hand, more recent studies of the MMOC involves MMOC with adjusted advection (MMOCAA) \cite{DFP97-MMOCAA-main}. Compared to the MMOC, MMOCAA has better global mass conservation properties, which is usually required for an accurate numerical simulation of models that are related to engineering problems. On the other hand, of particular difficulty in the implementation of ELLAM is an accurate evaluation of integrals involving steep back-tracked functions \cite{S15}. An inaccurate evaluation of these integrals will yield a loss in mass conservation; a fix to simplify the evaluation of these integrals, which will preserve the mass conservation property, has been proposed recently in \cite{AW11-stability-monotonicity-implementation}. 

Here, we present a detailed analysis in terms of the mass conservation properties, and provide an idea of combining the ELLAM and MMOC schemes, so that we can capitalise on the advantage of both schemes, while at the same time avoid their disadvantages. Having achieved global mass balance, a novel, less expensive adjustment yielding local volume conservation is then proposed. The complete coupled scheme then consists of a characteristic component (the combined ELLAM--MMOC), accompanied by a discretisation of the diffusive terms using the Gradient Discretisation Method (GDM) framework \cite{gdm}. This framework contains many classical methods (finite elements, finite volumes, etc.) for diffusion equations.
The complete coupled scheme, named GEM (for GDM--ELLAM--MMOC), therefore presents in one form several possible discretisations of the advection--diffusion model.

The main contributions of this work are
\begin{itemize}
	\item  a new adjustment algorithm to achieve local volume conservation (Section \ref{sec:local-mass-bal}),
	\item  precise analysis of mass balance errors for MMOC scheme (Section \ref{subsec mass-bal-MMOC}),
	\item  combined ELLAM--MMOC scheme for the advection-reaction equation~\eqref{eq:advection2} (Section \ref{sec:ELLAM-MMOC}),
	\item  application of the ELLAM--MMOC scheme and local volume adjustments for the miscible flow model (GEM scheme, Section \ref{sec:Peaceman}).
\end{itemize}  

The paper is made up of two major components. The first component focuses on the  advection-reaction equation \eqref{eq:advection2} and the characteristic based schemes used to discretise this equation. We start here by studying the ELLAM scheme and the MMOC scheme in Sections \ref{sec: ELLAM} and \ref{sec:MMOC} respectively. In particular, we look into the global mass balance properties of both schemes, and then discuss, for the ELLAM, how to achieve local volume conservation after achieving global mass balance. We then propose in Section \ref{sec:ELLAM-MMOC} a combined ELLAM--MMOC scheme, and discuss its advantages over both the ELLAM and the MMOC. The second component then focuses on the application of the proposed ELLAM--MMOC combination to a miscible flow model \eqref{eq:model} in porous medium. This model involves diffusion terms, that are discretised using the generic GDM framework. The local volume adjustments performed in Section \ref{sec:local-mass-bal} are then adapted for the numerical schemes for model \eqref{eq:model}. Numerical results presented in Section \ref{sec: num-results} illustrate the clear advantages of this GEM scheme. Finally, we provide in Section \ref{sec:conv.anal} a convergence result for the MMOC scheme. This can then be extended, together with the convergence result for the ELLAM scheme \cite{CDL17-convergence-ELLAM}, to obtain the convergence of the ELLAM--MMOC scheme for~\eqref{eq:advection}.

\subsection{Models}

Our objective is to design a robust, characteristic-based numerical scheme for a model of miscible displacement in a porous medium. This model, described in Section \ref{sec:Peaceman}, involves an elliptic equation for the pressure, and an advection--reaction--diffusion equation for the concentration of the invading fluid. For simplicity, we describe the characteristic-based scheme for the concentration equation without explicitly referring to the pressure equation. We therefore consider the scalar model
\begin{equation}\label{eq:advection}
\left\{\begin{array}{ll}
\phi \dfrac{\partial c}{\partial t} 
+ \div(\Daru c -\Lambda\nabla c) 
= f(c) &\qquad \mbox{ on } Q_T:=\O \times (0,T)\\
\Lambda\nabla c\cdot\bfn=0&\qquad\mbox{ on }\partial\O\times (0,T),\\
c(\cdot,0)=c_{\rm ini}&\qquad\mbox{ on }\Omega,
\end{array}\right.
\end{equation}
in which $T>0$, $\O$ is an open bounded domain of $\R^d$ ($d\ge 1$), the porosity $\phi$, the diffusion tensor $\Lambda$ and the velocity $\Daru$ are given, $\Daru\cdot\bfn=0$ on $\partial\O$, and $f(c)=f(c,\x,t)$ is a function $\R\times Q_T\to \R$. The unknown $c(\x,t)$ represents the amount of material (a fraction) present at $(\x,t)$. The characteristic method only deals with the advective part of the model, and will therefore be described on the advection--reaction equation (corresponding to $\Lambda \equiv 0$):
\begin{equation}\label{eq:advection2}
\left\{\begin{array}{ll}
\phi \dfrac{\partial c}{\partial t} 
+ \div(\Daru c) 
= f(c) &\qquad \mbox{ on } Q_T:=\O \times (0,T)\\
c(\cdot,0)=c_{\rm ini}&\qquad\mbox{ on }\Omega.
\end{array}\right.
\end{equation}
Note that the boundary is non-characteristic due to the assumption $\Daru\cdot\bfn=0$ on $\partial\O$, and thus no boundary conditions need to be enforced in \eqref{eq:advection2}.


\subsection{Assumptions on the data, and numerical setting} \label{sec:data-and-num-setting}

Throughout the article we assume the following properties:
\begin{subequations}\label{assump.global}
		\begin{equation} \label{assum:init-source}
		\begin{aligned}
			&c_{\rm ini}\in L^\infty(\O), \\
			&f:\R\times Q_T\to\R\mbox{ is Lipschitz continuous w.r.t.\ its first variable}\\
			&\mbox{and $f(0,\cdot,\cdot)\in L^\infty(Q_T)$},
		\end{aligned}
	\end{equation}
	\be
	\begin{aligned}
		&\phi\in L^\infty(\O)\mbox{ and there exists $\phi_*>0$ s.t. } \phi\ge \phi_*\mbox{ a.e.~on $\O$}\label{hyp:phi}
	\end{aligned}
	\ee
	\be\label{hyp:velocity}
	\begin{aligned}
	&\Daru \in L^\infty(0,T;L^2(\O)^d) \mbox{ and }\div \Daru\in L^\infty(Q_T).
	\end{aligned}
	\ee
\end{subequations}


Our objective in this paper is to describe numerical methods in a general setting, to ensure that our design and analysis of ELLAM--MMOC schemes applies at once to various possible spatial discretisations (e.g.\ finite-element or finite-volume based). To achieve this, we use the Gradient Discretisation Method (GDM), a generic framework for numerical methods for diffusion equations \cite{gdm}. Although most of our work will be done here on the advective--reactive parts of \eqref{eq:advection}, we will demonstrate that the GDM also provides all the required tools to describe ELLAM and MMOC schemes.

The GDM consists in replacing, in weak formulations of the models, the continuous (infinite-dimensional) spaces and corresponding operators by a discrete (finite-dimensional) space and reconstructions of functions and gradients. A space-time Gradient Discretisation (GD) is
$\discC=(X_\discC,\Pi_\discC,\nabla_\discC,\ICinterp_{\discC},(t^{(n)})_{n=0,\dots, N})$, where
\begin{itemize}
\item $X_{\discC}$ is a finite-dimensional real space, describing the unknowns of the
chosen scheme,
\item $\Pi_\discC:X_\discC\to L^\infty(\Omega)$ is a linear operator that reconstructs a function
from the unknowns,
\item $\nabla_\discC:X_\discC\to L^\infty(\Omega)^d$ is a linear operator that reconstructs a gradient
from the unknowns,
\item $\ICinterp_\discC c_{\rm ini}$ is a rule to interpolate $c_{\rm ini}$ onto an element
of $X_{\discC}$,
\item $(t^{(n)})_{n=0,\dots, N}$ are the time steps, and we let
$\dt^{(n+\frac12)}=t^{(n+1)}-t^{(n)}$.
\end{itemize}
Different choices of $\discC$ lead to different schemes, e.g.\ finite elements or finite volumes \cite{gdm}.

Finally, we assume in the following that $\Daru$ is approximated on each time interval $(t^{(n)},t^{(n+1)})$ by a function 
\begin{equation}\label{hyp:un}
\Daru^{(n+1)}\in L^2(\O)^d\mbox{ such that }\div\Daru^{(n+1)}\in L^\infty(\O).
\end{equation}
Precise assumptions regarding this approximation will be described in Section \ref{sec:conv.anal}.

In the rest of the paper, the variables are only made explicit in the integrals when there is a risk of confusion. Otherwise we simply write, e.g., $\int_\O qd\x$.

\section{ELLAM scheme for the advection--reaction equation} \label{sec: ELLAM}

\subsection{Motivation}
For any sufficiently smooth function $\varphi$, the product rule yields
\[
\varphi  \dfrac{\partial c}{\partial t}  =  \dfrac{\partial (c\varphi)}{\partial t} - c\dfrac{\partial \varphi}{\partial t}.
\]
Hence, \eqref{eq:advection2} gives, for any time interval $(t^{(n)},t^{(n+1)})$,
\begin{align*}
\int_{t^{(n)}}^{t^{(n+1)}}&\int_\O\phi(\x)\dfrac{\partial (c\varphi)}{\partial t} (\x,t) d\x dt\nonumber\\
&-
\int_{t^{(n)}}^{t^{(n+1)}}\int_\O c(\x,t)
\bigg[
\phi(\x)\dfrac{\partial \varphi}{\partial t} (\x,t) + \Daru(\x,t)\cdot \nabla \varphi(\x,t)
\bigg] d\x dt\nonumber\\
&=
\int_{t^{(n)}}^{t^{(n+1)}}\int_{\Omega} f(c,\x,t)\varphi(\x,t)d\x dt.
\end{align*}
To simplify the second term on the left hand side of the above equation,
the ELLAM requires that test functions $\varphi$ satisfy
\begin{equation}\label{eq:testfunc}
\phi\dfrac{\partial \varphi}{\partial t} + \Daru\cdot \nabla \varphi= 0 
\quad\text{ on }\O \times(t^{(n)},t^{(n+1)}),
\end{equation}
with $\varphi(\cdot,t^{(n+1)})$ given.
The advection equation~\eqref{eq:advection2} then leads to the relation
\begin{multline*}
\int_\O\phi(\x)(c\varphi) (\x,t^{(n+1)}) d\x
-
\int_\O\phi(\x)(c\varphi) (\x,t^{(n)}) d\x\\
=
\int_{t^{(n)}}^{t^{(n+1)}}\int_{\Omega} f(c,\x,t)\varphi(\x,t)d\x dt.
\end{multline*}

\subsection{ELLAM scheme}

The ELLAM scheme consists in exploiting the motivation above, in the discrete context of the GDM
in which trial and test functions are replaced by reconstructions $\Pi_\discC$ applied to trial and
test vectors in $X_\discC$.

\begin{definition}[ELLAM scheme]\label{def:ELLAM}
	Given a Gradient Discretisation $\discC$ and using a weighted trapezoid rule with weight $w\in [0,1]$ for the time-integration of the source term, the ELLAM scheme for \eqref{eq:advection2} reads as: find 
	 $(c^{(n)})_{n=0,\ldots,N}\in X_{\discC}^{N+1}$ such that $c^{(0)}=\ICinterp_{\discC} c_{\rm ini}$
	and, for all $n=0,\ldots,N-1$, 
		$c^{(n+1)}$ satisfies
		\begin{align}\label{conc-ELLAM}
		&\int_{\O} \phi \PiDc c^{(n+1)} \PiDc z -\int_{\O} \phi \PiDc c^{(n)} v_z(t^{(n)})\nonumber \\
		&\qquad= w\dtDisc \int_{\O} f_n v_z(t^{(n)}) +(1-w) \dtDisc \int_{\O} f_{n+1}\PiDc z
\qquad \forall z \in X_{\discC}, 
		\end{align}
		where $v_z$ is the solution to 
	\begin{equation}\label{Advection-test}
	\phi \partial_t v_z+\Daru^{(n+1)} \cdot \nabla v_z =0 \quad \text{ on } (t^{(n)},t^{(n+1)})\,,\mbox{ with $v_z(\cdot,t^{(n+1)})=\PiDc z$ }.
	\end{equation} 	
		Here and in the rest of the paper, we let $f_k:=f(\Pi_{\discC} c^{(k)},\cdot,t^{(k)})$ (or with
a suitable average over $(t^{(k)},t^{(k+1)})$ if $f$ is not continuous in time).
\end{definition}
\begin{remark}[About the time integration]
   The velocity field $\Daru$ was approximated by its value at time $t^{(n+1)}$, given by $\Daru^{(n+1)}$. Other choices for the approximation of $\Daru$, such as a centred approximation $\frac{1}{2}(\Daru^{(n)}+\Daru^{(n+1)})$, may be made, but we noticed in our tests that this does not noticeably change the numerical solution. A weighted trapezoid rule is applied for time integration in Definition \ref{def:ELLAM} for the purpose of achieving mass conservation, as discussed in \cite{AH06}.
	\end{remark}
Define the flow $F_{t}:\O\to\O$ such that, for a.e.\ $\x\in\O$,
\begin{equation} \label{charac}
\dfrac{dF_{t}(\x)}{dt} =\dfrac{\Daru^{(n+1)}(F_{t}(\x))}{\phi(F_t(\x))}\quad\mbox{ for $t\in [-T,T]$}, \qquad F_{0}(\x)= \x.
\end{equation}
Under Assumption \eqref{hyp:un}, the existence of this flow is proved in \cite[Lemma 5.1]{CDL17-convergence-ELLAM}.
The solution to \eqref{Advection-test} is then understood in the sense:
for $t\in (t^{(n)},t^{(n+1)}]$ and a.e.\ $\x\in\O$, $v_z(\x,t)=\Pi_\discC z(F_{t^{(n+1)}-t}(\x))$. In particular,
\begin{equation}\label{eq:testfunc.ELLAM}
v_z(\cdot,t^{(n)})=\PiDc z(F_{\dt^{(n+\frac12)}}(\cdot)).
\end{equation}

For any functions $f$ and $g$, defining the vector functions $f^{(n,w)}$ and $g_{F}$ by
\begin{equation}\label{def:fnw.gF}
\begin{aligned}
f^{(n,w)}(\x) &:= \bigg(wf(\x,t^{(n)}),(1-w)f(\x,t^{(n+1)})\bigg),\\
g_{F}(\x)&:= \bigg(g(F_{\dt^{(n+1/2)}}(\x)),g(\x)\bigg),
\end{aligned}
\end{equation}
the time-stepping~\eqref{conc-ELLAM} can be rewritten in the condensed form
\begin{equation}\label{ELLAM:condensed}
\int_{\O} \phi \PiDc c^{(n+1)} \PiDc z -\int_{\O} \phi \PiDc c^{(n)} v_z(t^{(n)})
=
\dtDisc\int_\O f^{(n,w)}\cdot (\PiDc z)_{F}.
\end{equation}

\subsection{Physical Interpretation}
We provide a simple physical interpretation of the ELLAM, by supposing to simplify that $\Pi_\discC$ is a piecewise-constant reconstruction on a given mesh $\mesh$. We also assume that for each cell $K\in\mesh$, there is $z_K\in X_\discC$ such that $\Pi_\discC z_K=\mathbbm{1}_K$. Writing $\Pi_{\discC} c^{(k)}= \sum_{M\in\mesh} c_M^{(k)} \mathbbm{1}_M$ and taking $z_K$ as test function, \eqref{ELLAM:condensed} and \eqref{eq:testfunc.ELLAM} give
\[
\begin{aligned}
\int_{K} \phi \PiDc c^{(n+1)} d\x={}&\int_{\O} \phi \sum_{M\in\mesh} c_M^{(n)} \mathbbm{1}_M(\x) \mathbbm{1}_K(F_{\dtDisc}(\x))d\x\\
&+\dtDisc\int_\O f^{(n,w)}\cdot (\mathbbm{1}_K)_{F}d\x,
\end{aligned}
\]
which reduces to
\begin{equation}\label{eq:phys-inter-ELLAM}
|K|_\phi c_K^{(n+1)}= \sum_{M\in\mesh}|M\cap F_{-\dtDisc}(K) |_{\phi}c_M^{(n)}+
\dtDisc\int_\O f^{(n,w)}\cdot (\mathbbm{1}_K)_{F}d\x,
\end{equation}
where $|E|_{\phi}=\int_E \phi$ is the available porous volume in a set $E\subset \R^d$.
The first term on the right hand side of \eqref{eq:phys-inter-ELLAM} tells us that the amount of material $c_K^{(n+1)}$ present in a particular cell $K\in\mesh$ at time $t^{(n+1)}$ is obtained by locating where the material in cell $K$ comes from, hence back-tracking the cell $K$ to $F_{-\dt^{(n+1/2)}}(K)$, measuring how much of the material $c^{(n)}_M$ is taken from each $M\in\mesh$, and deposing this material into the cell $K$. These are accompanied by the contribution of the source term $f$ in the particular cell $K$, which is given by the second term. We note here that this second term has a very similar treatment as the first term, i.e. the contribution that comes from $f$ at time $t^{(n)}$ is determined by the traceback region associated to cell $K$.

\subsection{Mass balance properties}
One desirable property for numerical schemes is conservation of mass. Essentially, we want a discrete form of the following equation, obtained by integrating \eqref{eq:advection2} over $\O$ and which tells us that  the change in $c$ is dictated by the amount of inflow/outflow given by the source term.
\begin{equation} \label{mass-bal}
\int_\O\phi(\x) c(\x,t^{(n+1)})d\x =\int_\O \phi(\x)c(\x,t^{(n)})d\x +\int_{t^{(n)}}^{t^{(n+1)}}\int_{\Omega} f(c,\x,t) d\x dt.
\end{equation}

Hence, it is important to define a measure of the mass balance error. Setting $\mathbf{e} :=(1,1)$,
the (discrete) mass balance error is defined by
	\begin{equation}\label{def:mass-bal-error}
	e_{\rm mass}:=\left|\int_\O\phi \PiDc c^{(n+1)}d\x
	-\int_\O \phi \PiDc c^{(n)}d\x 
	-\dtDisc\int_{\Omega} f^{(n,w)}\cdot \mathbf{e} \, d\x\right|.
	\end{equation}

A mass balance-preserving method is one for which $e_{\rm mass}=0$.
The ELLAM scheme \eqref{conc-ELLAM} satisfies this property. Indeed, taking $z_1=\sum_{K\in\mesh}z_K$, which satisfies $\Pi_\discC z_1=1$ over $\Omega$, as a test function in \eqref{ELLAM:condensed} gives $e_{\rm mass}=0$.


\begin{remark}[Steep back-tracked functions]\label{rem:ELLAM.steep}
The natural physical interpretation of ELLAM, together with its mass conservation property, seem to indicate that the ELLAM scheme should be preferred over other numerical schemes for the advection equation \eqref{eq:advection2}. However, for Darcy velocities typically encountered in reservoir simulations, the streamlines of the flow $F_t$ concentrate around injection wells, and the functions $v_z$ defined by \eqref{eq:testfunc.ELLAM} are then extremely steep in these regions. An accurate approximation of the integral of these functions in cells close to the injection well then requires to track a lot of quadrature points, which is very costly \cite{S15}. In some instances, even tracking several points along these regions would not give an accurate depiction of the integral. This is one of the main issues with ELLAM implementations. Fixes have been proposed, but they consist is resorting to a different approach, near the injection wells, than the ELLAM process \cite{AW11-stability-monotonicity-implementation}. We aim at designing a numerical scheme that readily behaves well, without having to implement specific fixes in certain regions. The MMOC will be instrumental to that objective.
\end{remark}

\subsection{Local volume conservation}\label{sec:local-mass-bal}

\begin{figure}[h]
	\centering
	\begin{tabular}{c@{\hspace*{2em}}c}
		\includegraphics[width=0.4\textwidth]{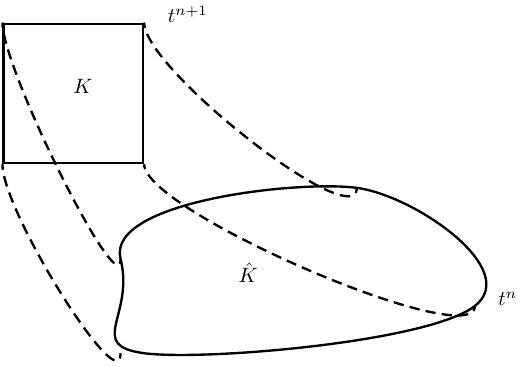} & 		\includegraphics[width=0.4\textwidth]{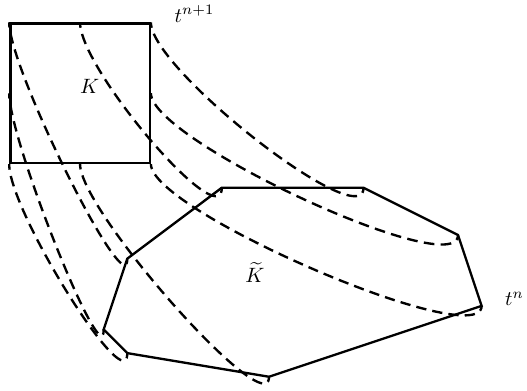}\\
	\end{tabular}
	\caption{ Traceback region $\widehat{K}$ (left: exact; right: approximation).}
	\label{trace-back-regions}
\end{figure}

The main difficulty of implementing an ELLAM type scheme is the evaluation of the integral $\int_{\widehat{K}} \phi \PiDc c^{(n)} d\x$ for each cell $K$, where $\widehat{K}=F_{-\dtDisc}(K)$. In general, the region $\widehat{K}$ (see Figure \ref{trace-back-regions}, left) cannot be exactly described, and hence we approximate it by polygons obtained from back-tracking the vertices, together with a number of points along the edges of the cell $K$. Figure \ref{trace-back-regions} (right) gives an illustration of the approximate traceback region $\widetilde{K}$ obtained by tracking the vertices, together with the edge midpoints of the cell $K$. 

In general, $|\widetilde{K}| \neq |F_{-\dtDisc}(K)|$. However, the equality of these volumes is essential in numerical simulations to preserve the accuracy of numerical solutions. To illustrate this point, consider the simple case of a divergence free velocity field in \eqref{eq:advection2}, with $\phi=1, f=0$ and $c_{\rm{ini}} = 1$. In this test case, the exact solution is given by $c(\x,t)=1$. In theory, upon implementing a numerical scheme with piecewise constant approximations for the unknown $c$, we should have the following simplified form of \eqref{eq:phys-inter-ELLAM} at the first time step:
\begin{equation}\nonumber
|K| c_K^{(1)} = \sum_{M\in\mesh}|M\cap F_{-\dtDisc}(K) |(c_{\rm{ini}})_M.
\end{equation}
However, due to the approximation of the traceback region, we only have 
\[
|K| c_K^{(1)} = \sum_{M\in\mesh}|M\cap \widetilde{K} |(c_{\rm{ini}})_M = |\widetilde{K}|,
\]
and thus 
\[
c_K^{(1)} = \dfrac{|\widetilde{K}|}{|K|}\neq 1
\]
in general. This example shows that an inaccurate approximation of the volume of the tracked cell renders the numerical scheme unable to recover constant solutions. Hence, we need to perform some adjustments on the polygonal region $\widetilde{K}$ in order to yield $|\widetilde{K}| = |F_{-\dtDisc}(K)|$, which we shall define as the \textit{local volume constraint} for $K$. In \cite{AH06,AH10-Fully-Conservative-ELLAM}, the local volume constraint is satisfied by moving through the tracked cells one by one, and adjusting the locations of the tracked midpoints along the direction of the characteristics. This was illustrated to work for square cells; however, there is no guarantee that such adjustments will yield a valid mesh configuration for irregular cells. A more generic approach is taken in \cite{D16-Opti-meshCorr}, where the local volume constraint is satisfied by solving an optimisation problem of finding a mesh which is closest to the tracked mesh formed by $\bigcup_{K\in\mesh} \widetilde{K}$, satisfying the local volume constraints. Common to these algorithms is an explicit expression for the adjusted traceback regions. However, as can be seen in \eqref{eq:phys-inter-ELLAM}, this is not necessary for piecewise constant approximations, standard in reservoir simulations based on finite volume methods. The important aspect is the computation of the quantities $|M\cap \widetilde{K}|$. 

\begin{figure}[h]
	\centering
	\begin{tabular}{c@{\hspace*{2em}}c}
		\includegraphics[width=0.4\textwidth]{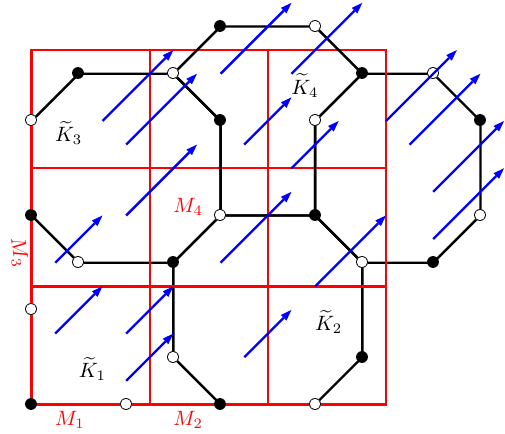} & 		\includegraphics[width=0.4\textwidth]{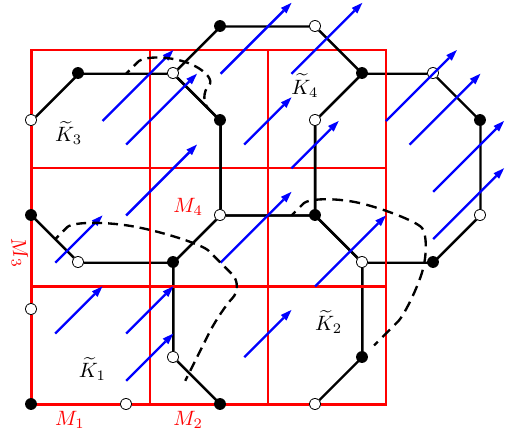}\\
	\end{tabular}
	\caption{ Traceback regions $\widetilde{K}_i$ (left: initial; right: illustration of possible perturbed cells after proposed volume adjustment).}
	\label{traceback-regions-adjustments}
\end{figure}

We propose an algorithm which adjusts $|M\cap \widetilde{K}|$ for each cell $K$, in the sense that these adjusted volumes would be something we expect to recover from a mesh obtained by adjusting the tracked points. The proposed algorithm works on any type of cells, but for simplicity of exposure we illustrate it in Figure \ref{traceback-regions-adjustments} using square cells, and traceback regions $\widetilde{K}_i$ approximated by tracking the vertices and edge midpoints of $K_i$. In Figure \ref{traceback-regions-adjustments}, the blue lines denote the trajectory taken by the velocity field, and the red squares form part of the original mesh cells $K\in\mesh$, whereas the black cells are their traceback regions: shaded points are the tracked vertices, and the hollow points are the tracked edge midpoints. In practice, after performing the tracking, aside from the final location of these vertices and midpoints, we also store the cell that they finally reside in. The algorithm is implemented cell-wise, starting from the first cell $K_1$, and proceeds as follows: Consider a cell $K_1$ with neighbors $K_2$, $K_3$, etc. This leads to a tracked cell $\widetilde{K}_1$ with neighbors $\widetilde{K}_2$, $\widetilde{K}_3$, etc. Suppose that $\widetilde{K}_1$ intersects the residing cells $M_1,M_2,M_3$ and $M_4$ (see Figure \ref{traceback-regions-adjustments}, left).  
\begin{itemize}
	\item[i)] We start by measuring the error $e_{K_1}:= |\widetilde{K}_1| - |F_{-\dtDisc}(K_1)|$. The relation $e_{K_1}>0$ (resp. $e_{K_1}<0$) means that we overestimate (resp. underestimate) the volume of the tracked region.
	\item[ii)]  We then compute the magnitude $|\Daru|$ of $\Daru$ at the tracked midpoints and also check whether $\Daru$ points into $\widetilde{K}_1$ or not. If the velocity points towards the same direction for two consecutive midpoints, then we also compute the magnitude of $\Daru$ for the vertex in between them.
	\item[iii)] We now illustrate how to adjust the volumes of the regions. If $e_{K_1}<0$, then it means that the volume $|\widetilde{K}_1|$ should be increased. Based on the velocity field $\Daru$ in Figure \ref{traceback-regions-adjustments}, this should be done by increasing along the direction of $\widetilde{K}_2$ and $\widetilde{K}_3$. Now, the velocity field along the edge midpoints that are located at $\widetilde{K}_1 \cap M_2$ and $\widetilde{K}_1 \cap M_3$ points outward of $\widetilde{K}_1$ and hence the vertex at $\widetilde{K}_1 \cap M_4$ should also be included. Denote then by $|\Daru_{1,i}|$ the magnitude of the velocity field evaluated at these tracked points in $\widetilde{K}_1 \cap M_i, i=2,3,4$. We will then adjust each of the volumes by subtracting, to each $|\widetilde{K}_1 \cap M_i|$, the quantity $\frac{|\Daru_{1,i}|}{\sum_{j=2}^{4}|\Daru_{1,j}|} e_{K_1}$. This will then make $K_1$ satisfy the local volume constraint.
	\item[iv)] To make sure that this quantity really represents something that would have come from a perturbed mesh (see Figure \ref{traceback-regions-adjustments}, right), the quantities $|\widetilde{K}_2 \cap M_2|, |\widetilde{K}_2 \cap M_4|, |\widetilde{K}_3 \cap M_3|$, and $|\widetilde{K}_3 \cap M_4|$ are adjusted accordingly (i.e. $\frac{|\Daru_{1,2}|}{\sum_{j=2}^{4}|\Daru_{1,j}|} e_{K_1}$ should be added onto $|\widetilde{K}_2 \cap M_2|$, and the corresponding quantities for the other edges). 
\end{itemize}  
We then proceed to adjust the other $\widetilde{K}_i$ in the same manner so that they satisfy their respective local volume constraints.

\begin{figure}[h]
	\centering
		\includegraphics[width=0.4\textwidth]{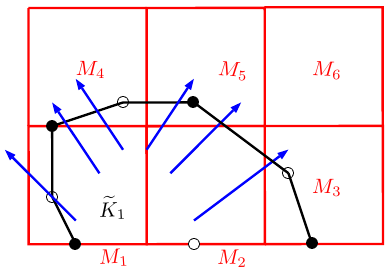}\\
	\caption{ Traceback region $\widetilde{K}_1$ of a cell at the boundary.}
	\label{traceback-regions-adjustments-bdry}
\end{figure}

As a remark, we note that one set of adjustments may not suffice to satisfy the local volume constraints for all cells. For example, if, after the adjustment of $\widetilde{K}_1$, we have that $e_{K_2}>0$, then we need to decrease the volume of $\widetilde{K}_2$. This can only be done (respecting the direction of the velocity field) by moving along the direction of $\widetilde{K}_1$. This will lead to $e_{K_1}>0$. Hence, after all the adjustments in the first stage, we need to re-evaluate $e_{K_i}$ and re-adjust the volumes. Of course, from a computational point of view, not all cells would be able to satisfy $e_{K_i}=0$ exactly, so we terminate our algorithm once $|e_{K_i}|$ is below a desired tolerance value for all cells. Another potential issue that may be encountered is when $K_1$ is a boundary cell (see Figure \ref{traceback-regions-adjustments-bdry}). If $\widetilde{K}_1$ lies on the boundary and $e_{K_1}>0$, then, we need adjustments which will decrease the volume of $\widetilde{K}_1$. Thus, adjusting along the direction of the velocity $\Daru$ will only worsen the problem, since it will further increase the volume of $\widetilde{K}_1$ and hence the value of $e_{K_1}$. In this case, the idea is to consider $-\Daru$ instead. This translates to pushing inwards the points that would have been pushed outwards if $e_{K_1}<0$. The re-adjustment of the other cells follow accordingly, and convergence is still expected. Making the volume adjustments in the direction of $\Daru$ (resp. $-\Daru$) corresponds to saying that we have tracked backward (resp. forward) too much, and hence to fix these, we must track forward (resp. backward) a little bit further.   

As has been mentioned in \cite{CD17-HMM-ELLAM-complete}, when dealing with irregular cells, we need to track more than just edge midpoints so that $\widetilde{K}$ is close to $F_{-\dtDisc}(K)$. Hence, for irregular cells, a slight modification for our algorithm should be made: we may opt to take smaller time steps (to make sure that the errors $e_{K_i}$ are small to start with, and adjustments can be made in a similar manner as with square cells, by only placing markers on the tracked midpoints and adjusting accordingly), or we may mark more than just the tracked midpoints (in order to have a better idea of the geometry of the tracked cell, and provide a more comprehensive list of the cell volumes to adjust). For our tests in Section \ref{sec: num-results}, when such a modification was required, we chose to take smaller time steps.

\section{MMOC scheme for the advection--reaction equation} \label{sec:MMOC}
\subsection{Motivation}
We use the product rule to write $\divg(\Daru c)=c\divg(\Daru)+\Daru \cdot \nabla c$. By treating $\phi\dfrac{\partial c}{\partial t} + \Daru \cdot \nabla c$ as a directional derivative, and denoting by $\tau$ the associated characteristic direction, we rewrite~\eqref{eq:advection2} as follows 
\begin{equation}\label{eq:concen1}
\zeta\dfrac{\partial c}{\partial \tau} 
=
f -c\,\divg(\Daru), 
\end{equation}  
where $\zeta= (\phi^2 + |\Daru|^2) ^{\frac{1}{2}}$.
The MMOC then considers the following approximation:
\begin{align}\label{eq:app deri}
\zeta(\x)\dfrac{\partial c}{\partial \tau} &\approx \zeta(\x) \dfrac{c(\x,t^{(n+1)})-c(\bar{\x},t^{(n)})}{((\x-\overline{\x})^2 + (\dtDisc)^2)^{\frac{1}{2}}}\nonumber\\
&\approx \phi(\x)\dfrac{c(\x,t^{(n+1)})-c(F_{-\dtDisc}(\x),t^{(n)})}{\dtDisc}.
\end{align} 
Here $\overline{\x}:= \x-\dfrac{\Daru^{(n)}(\x)}{\phi(\x)} \dtDisc $ and the flow $F_{t}$ is defined by \eqref{charac}. 

By integrating~\eqref{eq:concen1} over the time interval $[t^{(n)},t^{(n+1)}]$ and using the approximation~\eqref{eq:app deri} of the characteristic 
derivative, we obtain 
\begin{multline} \label{eq:app_adv}
\phi(\x)\bigg( c(\x,t^{(n+1)})-c(F_{-\dtDisc}(\x),t^{(n)}) \bigg)\\
\approx\int_{t^{(n)}}^{t^{(n+1)}} f(c(\x,t),\x,t)-c(\x,t) \divg \Daru^{(n+1)}(\x) dt.
\end{multline}
\subsection{MMOC scheme }

The MMOC scheme is written, in the GDM setting, by exploiting \eqref{eq:app_adv}.

\begin{definition}[MMOC scheme]
Given a Gradient Discretisation $\discC$ and using a weighted trapezoid rule with weight $w\in [0,1]$ for the time-integration of the source term, the MMOC scheme for \eqref{eq:advection2} reads as: find $(c^{(n)})_{n=0,\ldots,N}\in X_{\discC}^{N+1}$ such that $c^{(0)}=\ICinterp_{\discC} c_{\rm ini}$
	and, for all $n=0,\ldots,N-1$, $c^{(n+1)}$ satisfies
		\begin{align}\label{conc-MMOC}
		&\int_{\O} \phi \PiDc c^{(n+1)} \PiDc zd\x 
		-\int_{\O} \phi(\x) (\PiDc c^{(n)})(F_{-\dtDisc}(\x)) \PiDc z(\x) d\x\nonumber\\
		&\qquad= \dtDisc \int_{\O} \bigl[\bigl(f^{(n,w)}-(\Pi_{\discC} c)^{(n,w)} \divg \Daru^{(n+1)}\bigr)\cdot \mathbf{e}\bigr]\PiDc z d\x\qquad \forall z \in X_{\discC},
		\end{align}
where we recall that $\mathbf{e} :=(1,1)$, and where we have set (by slight abuse of the notation \eqref{def:fnw.gF})
\begin{equation}\label{def:Pic.nw}
(\Pi_\discC c)^{(n,w)}(\x) := \left(w\Pi_\discC c^{(n)}(\x),(1-w)\Pi_\discC c^{(n+1)}(\x)\right).
\end{equation}
\end{definition}
\subsection{Physical Interpretation}
As with the ELLAM, an interpretation of the MMOC will be provided for the simple case wherein we have a piecewise constant approximation of $c$. Fixing $K\in\mesh$ and taking in \eqref{conc-MMOC} the test vector $z_K$, such that $\Pi_{\discC} z_K= \mathbbm{1}_K$, we have
\[
\begin{aligned}
\int_{K} \phi {}&\PiDc c^{(n+1)} d\x =\int_{\O} \phi(\x) \sum_{M\in\mesh} c_M^{(n)} \mathbbm{1}_M(F_{-\dtDisc}(\x)) \mathbbm{1}_K(\x) d\x\\
&+
\dtDisc\int_K  f^{(n,w)}\cdot \mathbf{e} d\x - \dtDisc \int_{K} \bigl[\bigl(\Pi_{\discC} c\bigr)^{(n,w)} \divg \Daru^{(n+1)}\bigr]\cdot \mathbf{e} d\x,
\end{aligned}
\]
and thus
\begin{equation}\label{eq:phys-inter-MMOC}
\begin{aligned}
|K|_\phi {}&c_K^{(n+1)}= \sum_{M\in\mesh}|F_{\dtDisc}(M)\cap K |_\phi c_M^{(n)}\\
&+
\dtDisc\int_K f^{(n,w)}\cdot \mathbf{e} d\x - \dtDisc \int_{K} \bigl[(c_K )^{(n,w)}\divg \Daru^{(n+1)}\bigr]\cdot \mathbf{e} d\x.
\end{aligned}
\end{equation}
The first term on the right hand side of the equation tells us that the amount of material $c_K^{(n+1)}$ present in a particular cell $K\in\mesh$ at time $t^{(n+1)}$ is obtained by taking all cells $M\in\mesh$, advecting material from each of these cells (by computing the trace-forward regions $F_{\delta t ^{(n+1/2)}}(M)$), and determining which portion of each cell flows into $K$. The second term simply represents the change that comes from the source term $f$. We note that, unlike in the ELLAM, the contribution of the source term $f$ for the MMOC is taken exactly to be from cell $K$. By itself, this term tells us that, if the source term $f$ is nonconstant over regions close to one another, the MMOC will give either an excess or miss some amount that has flowed into the region $K$. The third term in \eqref{eq:phys-inter-MMOC} represents taking away a percentage of the net inflow/outflow in cell $K$ and, in some sense, attempts to balance out the excessive or missing amount resulting from the second term. 

\begin{remark}[Comparison of ELLAM and MMOC schemes]\label{rem:comp.ELLAM.MMOC}
	The ELLAM and the MMOC schemes are equivalent in a cell $K$ if the velocity field is divergence free, and the source term $f$ is constant in the region $F_{[-\dtDisc,0]}(K)=\cup_{t\in [-\dtDisc,0]}F_t(K)$. Physically, the equivalence is expected, as we are now just comparing the first terms of equations \eqref{eq:phys-inter-ELLAM} and \eqref{eq:phys-inter-MMOC}, which both compute the amount of substance that has flowed into cell $K$. This can be done in two ways: Either we first locate the regions from which the substance has come from (ELLAM), or we let the substances in all cells flow, and determine which ones enter the cell $K$ (MMOC). Mathematically, this equivalence can be established by performing a change of variables in $|F_{\dtDisc}(M)\cap K |_\phi$ and using the property $\phi(F_{\dtDisc}(\x))|JF_{\dtDisc}(\x)|=\phi(\x)$ if the flow occurs in a region where $\div\Daru^{(n+1)}=0$ (see \eqref{JacobianInt}).
\end{remark} 
\subsection{Analysis of mass balance error}\label{subsec mass-bal-MMOC}

Consider the MMOC scheme \eqref{conc-MMOC}. By taking the test function $z_1=\sum_{K\in\mesh}z_K$, we have $\Pi_\discC z_1=1$ in $\O$ and we obtain thus the discrete mass balance equation
\begin{equation}\label{mass-bal-MMOC}
\begin{aligned}
\int_\O  \phi\Pi_{\discC}c^{(n+1)}d\x 
={}& \int_{\O} \phi(\x) (\PiDc c^{(n)})(F_{-\dtDisc}(\x))  d\x\\
&+ \dtDisc \int_{\O} \bigl(f^{(n,w)}-(\Pi_{\discC} c)^{(n,w)} \divg \Daru^{(n+1)}\bigr)\cdot \mathbf{e}d\x.
\end{aligned}
\end{equation}
From this, we see that one of the disadvantages of MMOC schemes over ELLAM schemes is that, in general, MMOC schemes do not preserve mass. The mass balance error $e_{\rm mass}$ for MMOC is estimated by using~\eqref{mass-bal-MMOC} to substitute $\int_\O  \phi\Pi_{\discC}c^{(n+1)}d\x$ in \eqref{def:mass-bal-error}. Performing a change of variables, we obtain
\begin{align*}
e_{\rm mass}
={}&\bigg| \int_{\O} \phi(\x) (\PiDc c^{(n)})(F_{-\dtDisc}(\x))  d\x-\int_\O \phi \PiDc c^{(n)}d\x\\
&- \dtDisc \int_{\O} \bigl((\Pi_{\discC} c)^{(n,w)} \divg \Daru^{(n+1)}\bigr)\cdot \mathbf{e}d\x\bigg|\\
={}& \bigg| \int_{\O} \phi(F_{\dtDisc}(\x)) \PiDc c^{(n)}(\x) |JF_{\dtDisc}(\x)| d\x-\int_\O \phi \PiDc c^{(n)}d\x\\
&- \dtDisc \int_{\O} \bigl((\Pi_{\discC} c)^{(n,w)} \divg \Daru^{(n+1)}\bigr)\cdot \mathbf{e}d\x\bigg|.
\end{align*}
It follows from~\cite[Lemma 5.2]{CDL17-convergence-ELLAM} that
\begin{equation}\label{JacobianInt}
\phi(F_{\dtDisc}(\x))|JF_{\dtDisc}(\x)|- \phi(\x)
=\int_0^{\dtDisc} |JF_{t}(\x)| (\divg \Daru^{(n+1)})\circ F_{t}(\x) dt.
\end{equation}
Hence,
\begin{equation}
\label{mass-err-MMOC}
\begin{aligned}
e_{\rm mass}
={}& \bigg|\int_\O  \Pi_{\discC} c^{(n)}(\x) \int_0^{\dtDisc} |JF_{t}(\x)| (\divg \Daru^{(n+1)})\circ F_{t}(\x) dt d\x\\
&- \dtDisc \int_{\O} \bigl((\Pi_{\discC} c)^{(n,w)} \divg \Daru^{(n+1)}\bigr)\cdot \mathbf{e}d\x\bigg|\\
={}& \bigg| \int_0^{\dtDisc} \int_\O \Pi_{\discC}c^{(n)}(F_{-t}(\x)) \divg \Daru^{(n+1)}(\x)  d\x dt\\
&- \dtDisc \int_{\O} \bigl((\Pi_{\discC} c)^{(n,w)} \divg \Daru^{(n+1)}\bigr)\cdot \mathbf{e}d\x\bigg|.
\end{aligned}
\end{equation}
By the triangle inequality and recalling the definition \eqref{def:Pic.nw} of $(\Pi_{\discC} c)^{(n,w)}$, we infer 
\begin{equation}\label{MMOC:mass.balance}
\begin{aligned}
{}& e_{\rm mass}
\leq w\bigg| \int_0^{\dtDisc} \int_\O \bigg(\Pi_{\discC}c^{(n)}(F_{-t}(\x))-\Pi_{\discC} c^{(n)}(\x)\bigg) \divg \Daru^{(n+1)} (\x) d\x dt\bigg|\\
&+(1-w)\bigg| \int_0^{\dtDisc} \int_{\O} \bigl(\Pi_{\discC}c^{(n)}(F_{-t}(\x))-\Pi_{\discC} c^{(n+1)}(\x)\bigr) \divg \Daru^{(n+1)} (\x) d\x dt\bigg|. 
\end{aligned}
\end{equation}
This estimate shows that the mass balance error $e_{\rm mass}$ is minimal when $\dtDisc$ tends to $0$ (as $F_{-t}\to Id$ as $t\to 0$) or the approximate amount of substance $c$ 
near the non divergence-free regions, denoted by $U$, is almost constant. More precisely,
\[
\Pi_{\discC} c^{(n)} = \Pi_{\discC} c^{(n+1)} = \text{Const} \quad\text{on } F_{[-\dtDisc,0]}(U).
\]
Due to this, the MMOC is in most cases expected to give a relatively large mass balance error compared to ELLAM. Hence, it is not practical to perform the local adjustments described in Section \ref{sec:local-mass-bal} for the MMOC scheme.
\begin{remark}[Conservation of mass for the MMOC]\label{rem:massConsMMOC}
	Estimate \eqref{MMOC:mass.balance} shows that if the velocity field is divergence free then the MMOC scheme conserves mass, which is consistent with Remark \ref{rem:comp.ELLAM.MMOC}.
\end{remark}

\begin{remark}[Forward tracking and cost near the injection cells]
Contrary to the ELLAM, the MMOC requires to forward-track test functions (see \eqref{eq:phys-inter-MMOC} in the case of piecewise constant approximations). Hence, in the MMOC, functions whose support is near the injection wells are not backward tracked into the injection cells, which makes them very steep and difficult to integrate (see Remark \ref{rem:ELLAM.steep}), but forward tracked far from these cells into non-steep functions that are easier to integrate.
\end{remark}

\section{A combined ELLAM--MMOC scheme for the advection--reaction equation} \label{sec:ELLAM-MMOC}

Here, we propose a combined ELLAM--MMOC scheme, to benefit from the mass balance property of the ELLAM and mitigate its costly implementation near the injection wells by using the MMOC method, much less expensive in these regions.

We start by applying a pure ELLAM scheme over the first few time steps, until $c$ is almost constant in areas near the non divergence-free regions. After which, we do a split ELLAM--MMOC scheme, where we apply MMOC over these areas, and ELLAM elsewhere. The interest of such a scheme is twofold. First, the computational cost is reduced compared to a pure ELLAM scheme as we no longer have to compute integrals of steep functions. Second, upon using MMOC only in regions where $\div\Daru=0$ or $c$ is already almost constant, no mass balance error occurs. This combined scheme takes out the main disadvantages of both methods.

\subsection{Presentation of the ELLAM--MMOC scheme}

Take $\alpha$ a function of the space variable, write $c=\alpha c+(1-\alpha)c$ and decompose the model~\eqref{eq:advection2} into
\be\label{adv.2}
\phi\frac{\partial (\alpha c)}{\partial t} + \div((\alpha c)\Daru)
+\phi\frac{\partial ((1-\alpha) c)}{\partial t} + \div(((1-\alpha) c)\Daru)
=\alpha f + (1-\alpha)f.
\ee
Discretise this by applying ELLAM on the first part $\alpha c$ (and $\alpha f$) and MMOC on the second part $(1-\alpha)c$ (and $(1-\alpha)f$). Defining $\PiDc (\alpha c)^{(n)}$ as $\alpha \PiD c^{(n)}$, one time step of this leads to
\begin{align*}
&\int_{\O} \phi \PiDc c^{(n+1)} \PiDc zd\x -\int_{\O} \phi(\x) \alpha(\x)\PiDc c^{(n)}(\x) \PiDc z (F_{\dtDisc}(\x))d\x\\
&-\int_{\O} \phi(\x) \left[(1-\alpha)\PiDc c^{(n)}\right](F_{-\dtDisc}(\x)) \PiDc z (\x)d\x\\
&\quad= \dtDisc\int_\O \alpha f^{(n,w)}\cdot (\PiDc z)_Fd\x
+\dtDisc\int_\O \bigl[(1-\alpha) f^{(n,w)}\cdot \mathbf{e}\bigr]\PiDc zd\x\\
&\qquad
-\dtDisc \int_{\O} \bigl[(1-\alpha)\divg \Daru^{(n+1)}\bigl(\Pi_{\discC} c \bigr)^{(n,w)}\cdot \mathbf{e}\bigr] \PiDc zd\x.
\end{align*}

\begin{remark}[Interpretation of the combined ELLAM--MMOC]
An interpretation can be given by considering $c_1=\alpha c$ and $c_2=(1-\alpha)c$ as two miscible fluids (that are also miscible in their surroundings) that do not react with each other, and are advected by the velocity $\Daru$. We can consider the combination of these two as one single fluid with concentration $c$, that is advected at velocity $\Daru$ (one can also consider that $c_1$ is made of red molecules, $c_2$ of green molecules, in which case the combination $c$ is yellow; to advect this yellow fluid, one can advect the red molecules with $\Daru$ and the green ones with $\Daru$ too).  The presentations \eqref{eq:advection2} or \eqref{adv.2} correspond to one or the other of these interpretations: do we want to consider both fluids together, or do we treat them separately. For the numerical method, it consists in applying ELLAM on one and MMOC on the other.
\end{remark}

The following definition summarises the combined ELLAM--MMOC scheme.

\begin{definition}[ELLAM--MMOC scheme] \label{def:GEM.ar}
	Given a Gradient Discretisation $\discC$ and using a weighted trapezoid rule with weight $w\in [0,1]$ for the time-integration of the source term, the ELLAM--MMOC scheme for \eqref{eq:advection2} reads as: find $(c^{(n)})_{n=0,\ldots,N}\in X_{\discC}^{N+1}$ such that $c^{(0)}=\ICinterp_{\discC} c_{\rm ini}$
	and, for all $n=0,\ldots,N-1$, $c^{(n+1)}$ satisfies
	\begin{equation}
	\label{GSconc-ELLAM-MMOC}
	\begin{aligned}
	\int_{\O} \phi {}&\PiDc c^{(n+1)} \PiDc zd\x 
	-\int_{\O} \phi(\x) \alpha(\x)\PiDc c^{(n)}(\x) \PiDc z(F_{\dtDisc}(\x))  d\x\\
	 -\int_{\O} {}&\phi(\x)\left[(1-\alpha)\PiDc c^{(n)}\right](F_{-\dtDisc}(\x)) \PiDc z(\x) d\x\\
	={}& 
	\dtDisc \int_{\O} \alpha f^{(n,w)}\cdot (\PiDc z)_F
	+\dtDisc \int_{\O} \bigl[(1-\alpha) f^{(n,w)}\cdot \mathbf{e}\bigr] \PiDc z\\
	&-\dtDisc \int_{\O} \bigl[(1-\alpha)\divg \Daru^{(n+1)}\bigl(\Pi_{\discC} c \bigr)^{(n,w)}\cdot \mathbf{e}\bigr] \PiDc z\qquad\forall z \in X_{\discC}.
	\end{aligned}
	\end{equation}
\end{definition}
\subsection{Analysis of mass balance error}\label{subsec mass-bal-ELLAM-MMOC}
Taking $z_1=\sum_{K\in\mesh}z_K$ (so that $\PiDc z_1=1$ in $\O$) in \eqref{GSconc-ELLAM-MMOC} and plugging into \eqref{def:mass-bal-error}, the mass balance error $e_{\rm mass}$ of the ELLAM--MMOC scheme is estimated as follows:
\begin{align*}
e_{\rm mass}
={}&
\bigg|
\int_{\O} \phi(\x) \left[(1-\alpha)\PiDc c^{(n)}\right](F_{-\dtDisc}(\x))  d\x\\
&-\int_\O \phi (1-\alpha)\PiDc c^{(n)}d\x-\dtDisc \int_{\O} (1-\alpha) \divg \Daru^{(n+1)}\bigl(\Pi_{\discC} c \bigr)^{(n,w)}\cdot \mathbf{e}\,d\x\bigg|\\
={}&
\bigg| \int_{\O} \phi\left((F_{\dtDisc}(\x)\right) \left[(1-\alpha)\PiDc c^{(n)}\right](\x) |JF_{\dtDisc}(\x)|  d\x\\
&-\int_\O \phi (1-\alpha)\PiDc c^{(n)} d\x
-\dtDisc \int_{\O} (1-\alpha)\divg \Daru^{(n+1)}\bigl(\Pi_{\discC} c \bigr)^{(n,w)}\cdot \mathbf{e} \,d\x\bigg|.
\end{align*} 
By using~\eqref{JacobianInt} and doing a change of variable $F_{-t}$ as in \eqref{mass-err-MMOC}, we obtain
\begin{align*}
e_{\rm mass}
={}& \bigg| \int_0^{\dtDisc} \int_\O \left[(1-\alpha)\PiDc c^{(n)}\right](F_{-t}(\x)) \divg \Daru^{(n+1)}(\x)  d\x dt\nonumber\\
&-\dtDisc \int_{\O} (1-\alpha) \divg \Daru^{(n+1)}\bigl(\Pi_{\discC} c \bigr)^{(n,w)}\cdot \mathbf{e}\,d\x\bigg|.\\
\leq{}&
w \bigg|\int_0^{\dtDisc} \int_\O \bigg[ \bigl((1-\alpha) \PiDc c^{(n)}\bigr)(F_{-t}(\x)) \\
&\qquad\qquad\qquad\qquad\qquad- \bigl((1-\alpha) \PiDc c^{(n)}\bigr)(\x)\bigg] \divg \Daru^{(n+1)}(\x) d\x dt\bigg|\nonumber\\
&+(1-w)
\bigg|\int_0^{\dtDisc} \int_\O \bigg[ \bigl((1-\alpha) \PiDc c^{(n)}\bigr)(F_{-t}(\x)) \\
&\qquad\qquad\qquad\qquad\qquad- \bigl((1-\alpha) \PiDc c^{(n+1)}\bigr)(\x)\bigg] \divg \Daru^{(n+1)}(\x) d\x dt\bigg|.
\end{align*}
Hence, the mass balance error $e_{\rm mass}$ of the ELLAM--MMOC scheme is minimal when $\dtDisc$ tends to $0$ or, setting $U=\{\x\in\O\,:\,\div\Daru^{(n+1)}(\x)\not=0\}$, if
\begin{equation*} 
\bigl(1-\alpha \bigr) \Pi_{\discC} c^{(n)}
\approx 
\bigl(1- \alpha\bigr) \Pi_{\discC} c^{(n+1)}\approx \text{Const}\ 
\mbox{ on $F_{[-\dtDisc,0]}(U)$.}
\end{equation*}

\begin{remark}[mass conserving $\alpha$]\label{rem.ex:ELLAM-MMOC}
	In particular, mass conservation is achieved if 
	\begin{itemize}
		\item $\alpha= 1$ on $F_{[-\dtDisc,0]}(U)$ (that is, pure ELLAM is used on the traceback of non-divergence free regions), or
		\item $\Pi_{\discC} c^{(n)}\approx \Pi_{\discC} c^{(n+1)}\approx C_1$ and $\alpha\approx C_2$, where each $C_i$ is a constant, on 
		$$
		D:=\{ \x \in \O\,:\, \alpha(\x)\not= 1\}\cap F_{[-\dtDisc,0]}(U)
		$$
		(that is, if MMOC is used --partially or entirely-- on a domain $D$ that is inside the traceback of non-divergence free regions, then the approximate concentration should almost be constant and stationary on $D$, and $\alpha$ should also be constant on $D$).
	\end{itemize}
\end{remark}

\subsection{Implementation for piecewise constant test functions}

As with the ELLAM and MMOC, we consider a piecewise constant approximation for $c$. Then, considering the ELLAM--MMOC scheme in \eqref{GSconc-ELLAM-MMOC}, we write $\PiDc c^{(n)}=\sum_{M\in\mesh}c_M^{(n)} \mathbbm{1}_M$ and find $(c^{(n)})_{n=0,\ldots,N}\in X_{\discC}^{N+1}$ such that $c^{(0)}=\ICinterp_{\discC} c_{\rm ini}$
and, for all $n=0,\ldots,N-1$, $c^{(n+1)}$ satisfies
\begin{equation}\nonumber
\begin{aligned}
\int_{\O} \phi {}& c_{K}^{(n+1)}\mathbbm{1}_{K}d\x 
-\int_{\O} \phi(\x) \alpha(\x)\sum_{M\in\mesh}c_M^{(n)} \mathbbm{1}_M(\x) \mathbbm{1}_{K}(F_{\dtDisc}(\x))  d\x\\
-\int_{\O} {}&\phi(\x)\left(1-\alpha\right)(F_{-\dtDisc}(\x)) \sum_{M\in\mesh}c_M^{(n)} \mathbbm{1}_M(F_{-\dtDisc}(\x)) \mathbbm{1}_{K}(\x) d\x\\
={}& 
\dtDisc \int_{\O} \alpha f^{(n,w)}\cdot (\mathbbm{1}_{K})_F
+\dtDisc \int_{\O} \bigl[(1-\alpha)f^{(n,w)}\cdot \mathbf{e}\bigr] \mathbbm{1}_{K}\\
&-\dtDisc \int_{\O} \bigl[(1-\alpha)\divg \Daru^{(n+1)}\bigl(\Pi_{\discC} c \bigr)^{(n,w)}\cdot \mathbf{e}\bigr] \mathbbm{1}_{K} \quad \forall K\in\mesh.
\end{aligned}
\end{equation}
Assume that $\alpha$ is piecewise constant on $\mesh$ and only takes the values $0$ and $1$. Each cell $M\in\mesh$ can then be classified as $\mesh_{\rm ELLAM}$ (corresponding to $\alpha=1$) or $\mesh_{\rm MMOC}$ (corresponding to $\alpha=0$). The above relation is then re-written
\begin{equation}\label{GSconc-ELLAM-MMOC-pcwise-const}
\begin{aligned}
&c_{K}^{(n+1)}|K|_\phi \\
&\ - \sum_{M\in \mesh_{\rm ELLAM}} c_M^{(n)} |M\cap F_{-\dtDisc}(K)|_\phi- \sum_{M\in \mesh_{\rm MMOC}} c_M^{(n)} |F_{\dtDisc}(M)\cap K|_\phi\\
&\quad\quad =
\dtDisc \int_{\O} \alpha f^{(n,w)}\cdot (\mathbbm{1}_{K})_F
+\dtDisc \int_{\O} \bigl[(1-\alpha)f^{(n,w)}\cdot \mathbf{e}\bigr] \mathbbm{1}_{K}\\
&\quad\qquad -\dtDisc \int_{\O} \bigl[(1-\alpha)\divg \Daru^{(n+1)}\bigl(\Pi_{\discC} c \bigr)^{(n,w)}\cdot \mathbf{e}\bigr] \mathbbm{1}_{K} \quad \forall K\in\mesh.
\end{aligned}
\end{equation}
If $K\in \mesh_{\rm ELLAM}$, then we only need to compute the integral of the first term on the right hand side of \eqref{GSconc-ELLAM-MMOC-pcwise-const} since the latter terms will be zero. Otherwise, only the second and third terms are computed. These are approximated by taking the average value of $f$ and $\divg \Daru^{(n+1)}$ on the respective cells.


We will demonstrate in Section \ref{sec:Peaceman} that, with a proper choice of $\alpha$, the combined ELLAM--MMOC scheme can be implemented with an equivalent or cheaper computational cost
than ELLAM, does not degrade much (or at all) the global mass conservation properties (contrary to MMOC), and has a much better local mass conservation compared to ELLAM on non-Cartesian meshes.

\subsection{Comparison with the MMOCAA}
 Of particular interest is a comparison with the MMOC scheme with adjusted advection (MMOCAA), first introduced in \cite{DFP97-MMOCAA-main}. The MMOCAA is a modification of MMOC designed to conserve the global discrete mass. Simply stated, the modification consists of perturbing the foot of the characteristic $F_{-\dtDisc}(\x)$ by a term of order $O((\Delta t)^2)$. Fixing $\eta \in (0,1)$, set
	\begin{equation}\nonumber
	\begin{aligned}
	\bar{\x}^+ &=F_{-\dtDisc}(\x) + \eta \dfrac{\Daru^{(n+1)}(\x)}{\phi(\x)} (\dtDisc)^2 \\
	\bar{\x}^- &=F_{-\dtDisc}(\x) - \eta \dfrac{\Daru^{(n+1)}(\x)}{\phi(\x)} (\dtDisc)^2.
	\end{aligned}
	\end{equation}
	For simplicity of notation, we denote the difference in mass for the MMOC scheme  $d_{\rm mass}$ to be \[
		\begin{aligned} 
		d_{\rm mass}:=& \int_\O  \Pi_{\discC} c^{(n)}(\x) -\int_{\O} \phi(\x) (\PiDc c^{(n)})(F_{-\dtDisc}(\x)) \\
		&- \dtDisc \int_{\O} \bigl((\Pi_{\discC} c)^{(n,w)} \divg \Daru^{(n+1)}\bigr)\cdot \mathbf{e}d\x.
		\end{aligned}\] We then define

\[
		 \widehat{\Pi_\discC c^{(n)}}(\x):=\begin{cases}
	\max(\Pi_\discC c^{(n)}(\bar{\x}^+),\Pi_\discC c^{(n)}(\bar{\x}^-)) &\qquad \mbox{if $d_{\rm mass}\leq 0$ } \\
	\min(\Pi_\discC c^{(n)}(\bar{\x}^+),\Pi_\discC c^{(n)}(\bar{\x}^-)) &\qquad \mbox{otherwise.}
	\end{cases}
\]

	In order to enforce a discrete conservation of mass, the term $\Pi_\discC c(F_{-\dtDisc},t^{(n)})$ in \eqref{conc-MMOC} is then replaced by $\Pi_\discC c_{\gamma}(\x,t^{(n)}):=\gamma \Pi_\discC c(F_{-\dtDisc},t^{(n)})+(1-\gamma) \widehat{\Pi_\discC c^{(n)}}(\x)$, where $\gamma$ is chosen so that 
	$d_{\rm mass}$, with $\Pi_\discC c(F_{-\dtDisc},t^{(n)})$ replaced by $\Pi_\discC c_{\gamma}(\x,t^{(n)})$, is equal to 0. For a more detailed presentation and implementation of the MMOCAA, we refer the readers to \cite{DFP97-MMOCAA-main,DHP99-MMOCAA-analysis,H00-MMOCAA}.
	
	It should be noticed that, contrary to the underlying principles of the ELLAM--MMOC scheme, there is no physical justification for using this parameter $\gamma$ to enforce the mass conservation (that is, $d_{\rm mass}=0$). Moreover, in some instances, at the first few time steps of a simulation, $\int_{\O} \phi(\x) (\PiDc c^{(n)})(F_{-\dtDisc}(\x))\PiDc z(\x)  d\x=\int_{\O} \phi(\x)\widehat{\Pi_\discC c^{(n)}}(\x)\PiDc z(\x) d\x$ and thus mass conservation cannot be achieved for any $\gamma$; see \cite{DHP99-MMOCAA-analysis}. Also, in order to be able to determine the proper value for $\gamma$, one needs to evaluate both
\[
\int_{\O} \phi(\x) (\PiDc c^{(n)})(F_{-\dtDisc}(\x))\PiDc z(\x)  d\x\mbox{ and }\int_{\O} \phi(\x) \widehat{\Pi_\discC c^{(n)}}(\x)\PiDc z(\x)  d\x.
\]
However, even though $\gamma$ ensures that the MMOCAA achieves global mass balance, local mass conservation is not necessarily satisfied. Hence, after computing $\gamma$, adjustments, such as those at the core of our ELLAM--MMOC approach (see Section \ref{sec:local-mass-bal}), need to be performed in order for the MMOCAA to provide numerical results that preserve local mass balance.


\section{Application: miscible flow model} \label{sec:Peaceman}

\subsection{Presentation of the model}

We consider the model of miscible flow in porous medium represented by the following coupled system of elliptic and parabolic PDEs \cite{PR62,E83-Mathematics-Reservoir}:
\begin{subequations} \label{eq:model}
	\begin{equation} \label{eq:press}
	\begin{aligned}
	\nabla \cdot \textbf{u} &= q^{+}-q^{-} \qquad (\textbf{x},t) \in \Omega \times [0,T] \\
	\textbf{u} &= - \dfrac{\textbf{K}}{\mu(c)} \nabla p \qquad (\textbf{x},t) \in \Omega \times [0,T] \\
	\end{aligned}
	\end{equation}
	\begin{equation} \label{eq:conc}
	\phi \dfrac{\partial c}{\partial t} + \nabla \cdot (\textbf{u}c-\textbf{D}(\textbf{x},\textbf{u})\nabla c) = q^{+}-cq^{-} \qquad (\textbf{x},t) \in \Omega \times [0,T].
	\end{equation}
	The unknowns  are $p(\x,t)$, $\darcyU(\x,t),$ and $c(\x,t)$ which denote the pressure of the mixture, the Darcy velocity, and the concentration of the injected solvent, respectively.
	The functions $q^{+}$ and $q^{-}$ represent the injection and production wells respectively, and $\textbf{D}(\textbf{x},\Daru)$ is the diffusion--dispersion tensor
	\begin{equation} \nonumber
	\textbf{D}(\textbf{x},\Daru) = \phi(\textbf{x})\left[d_{m}\textbf{I}+d_{l}|\Daru|\proj(\Daru)+d_{t}|\Daru|\left(\textbf{I}-\proj(\Daru)\right)\right]
	\mbox{ with }\proj(\Daru) = \left(\dfrac{u_{i}u_{j}}{|\Daru|^{2}}\right)_{i,j}.
	\end{equation}
	Here, $d_{m}>0$ is the molecular diffusion coefficient, $d_{l}>0$ and $d_{t}>0$ are the longitudinal and transverse dispersion coefficients respectively, and $\proj(\Daru)$ is the projection matrix along the direction of $\Daru$.
	The absolute permeability is $\K$, a space-dependent symmetric, bounded uniformly coercive diffusion tensor, and $\mu(c)=\mu(0)[(1-c)+M^{1/4}c]^{-4}$ is the viscosity of the fluid mixture, where $M=\mu(0)/\mu(1)$ is the mobility ratio of the two fluids.
	As usually considered in numerical tests, we take no-flow boundary conditions:
	\begin{equation}
	\Daru \cdot \bfn = (\textbf{D}\nabla c) \cdot \bfn = 0 
	\mbox{ on }\partial\O \times (0,T).
	\end{equation}
\end{subequations}
The concentration equation is completed by an initial condition, and the pressure equation by
an average condition:
\[
c(\x,0)=c_{\rm ini} \mbox{ for all $\x\in \O$, }\quad\int_\O p(\x,t)d\x=0\mbox{ for all $t\in (0,T)$}.
\]

\subsection{GDM--ELLAM--MMOC (GEM) scheme}\label{sec: GEM scheme}
The idea is to implement a time-marching algorithm, wherein gradient discretisations (as described in Section \ref{sec:data-and-num-setting}) are used to approximate the diffusive terms for both \eqref{eq:press} and \eqref{eq:conc}. Some examples for which different GDs are applied for each equation in \eqref{eq:model} have been presented in \cite{CDL17-convergence-ELLAM}. Note that the GDs used for the pressure equation do not need to involve the time components (time steps and interpolant of the initial condition).
As highlighted in Section \ref{sec:ELLAM-MMOC}, combining the ELLAM and MMOC for the treatment of the advective terms removes the main disadvantages of each of the schemes. Hence, we propose to use the combined ELLAM--MMOC scheme for the advective component. We will refer to the combination of the GDM with the ELLAM--MMOC scheme for the complete coupled model \eqref{eq:model} as the GDM--ELLAM--MMOC (GEM) scheme. 

Writing $A=\dfrac{\mathbf{K}}{\mu}$, the following definition of the GEM scheme is inspired by the construction of the GDM--ELLAM scheme in \cite{CDL17-convergence-ELLAM,CD17-HMM-ELLAM-complete} and by the design of the ELLAM--MMOC scheme for the advection--reaction model (Definition \ref{def:GEM.ar}).

\begin{definition}[GEM scheme] \label{def:GEM}
		Let $\discP  = (X_{\discP}, \PiDp, \gradDp)$ be a space GD for the pressure, and $\discC=(X_\discC,\Pi_\discC,\nabla_\discC,\ICinterp_{\disc},(t^{(n)})_{n=0,\dots, N})$ be a time-space GD for the concentration. Let $\alpha:\O\to [0,1]$.
 The GEM scheme for \eqref{eq:model} reads as: find $(p^{(n)})_{n=1,\ldots,N}\in X_{\discP}^N$
	and $(c^{(n)})_{n=0,\ldots,N}\in X_{\discC}^{N+1}$ such that $c^{(0)}=\ICinterp_{\discC} c_{\rm ini}$
	and, for all $n=0,\ldots,N-1$,
	\begin{itemize}[leftmargin=1.5em]
		\item[i)] $p^{(n+1)}$ solves
		\begin{equation}\label{GSpress}
		\begin{aligned}
		& \int_{\O} \PiDp p^{(n+1)} =0 \mbox{ and }\\
		&\int_{\O} A(\x,\PiDc c^{(n)}) \gradDp p^{(n+1)} \cdot \gradDp z = \int_{\O}  (q^{+}_n -q^{-}_n) \PiDp z\,,\qquad
		\forall z \in X_{\discP}
		\end{aligned}
		\end{equation}
		where $q^\pm_n(\cdot)=\frac{1}{\dtDisc}\int_{t^{(n)}}^{t^{(n+1)}}q^\pm(\cdot,s)ds$
		(or, alternatively, $q^\pm_n=q^\pm(t^{(n)})$ if $q^\pm$ are continuous in time).
		
		\item[ii)] A Darcy velocity $\discDarcyU$ is reconstructed from $p^{(n+1)}$
		and, to account for the advection term in the concentration equation,
		the following advection equation is considered; it defines space-time
		test functions from chosen final values:
		\begin{equation}\label{Advection}
		\phi \partial_t v+\discDarcyU \cdot \nabla v =0 \quad \text{on } (t^{(n)},t^{(n+1)})\,,\mbox{ with $v(\cdot,t^{(n+1)})$ given}.
		\end{equation}
		
		\item[iii)] Using a weighted trapezoid rule with weight $w\in [0,1]$ for the time-integration of the source term and setting $\U_\discP^{(n+1)} = A(\x, \PiDc c^{(n)}) \gradDp p^{(n+1)}$, $c^{(n+1)}$ satisfies
	\begin{align}\label{GSconc-GEM}
\int_{\O} \phi{}&\PiDc c^{(n+1)} \PiDc zd\x 
-\int_{\O} \phi(\x) \left[\alpha\PiDc c^{(n)}\right](\x) \PiDc z(F_{\dtDisc}(\x)) d\x\nonumber\\
-\int_{\O}& \phi(\x) \left[(1-\alpha)(\PiDc c^{(n)})\right](F_{-\dtDisc}(\x)) \PiDc z(\x) d\x\nonumber\\
&+\dtDisc\int_\O \mathbf{D}(\mathbf{U}_\discP^{(n+1)})\nabla_\discC c^{(n+1)}\cdot\nabla_\discC z\nonumber\\
={}& 
\dtDisc \int_{\O} \alpha \left[(q^+ -\PiDc c q^-)\right]^{(n,w)}\cdot (\PiDc z)_F\nonumber\\
&+\dtDisc \int_{\O} \bigl[(1-\alpha)\bigl(q^+ (1-\PiDc c)\bigr)^{(n,w)}\cdot \mathbf{e}\bigr] \PiDc z,
\qquad\forall z \in X_{\discC},
\end{align}
where we recall the notations \eqref{def:fnw.gF}, and we set $q^\pm_N=q^\pm_{N-1}$ if these quantities are defined by averages on time intervals (there is no time interval $(t^{(N)},t^{(N+1)})$).
	\end{itemize}
\end{definition} 

\begin{remark}[About the gradient discretisation and the velocity reconstruction]
	The gradient discretisation $\discP$ for the pressure equation \eqref{GSpress} can be chosen so that the scheme is locally conservative. Some examples of such schemes are the mixed finite element (MFEM) and hybrid mimetic mixed (HMM) schemes. If the gradient discretisation gives a velocity field $A(\x,\PiDc c^{(n)}) \gradDp p^{(n+1)}$ that is already in $H(\divg,\O)$, which is the case if the GD corresponds to an MFEM scheme, then the Darcy velocity $\discDarcyU$ is taken equal to that field. For some other gradient discretisations, such as the ones corresponding to the HMM scheme, the velocity field is not in  $H(\divg,\O)$ and a specific reconstruction, based on the numerical fluxes, must be made (see e.g. \cite[Section 2.3]{CD17-HMM-ELLAM-complete}).
	\end{remark}

Key to an efficient and accurate implementation of the GEM scheme is a proper choice of $\alpha$ so that mass conservation is achieved without having to deal with the steep source terms encountered in ELLAM. Remarks \ref{rem:ELLAM.steep} and \ref{rem.ex:ELLAM-MMOC} give us an idea of how to define the function $\alpha$. In the context of the complete coupled model \eqref{eq:model}, the non-divergence free regions are the injection and production cells. Moreover, it is expected that, for an injection cell $C_+$, $F_{[-\dtDisc,0]}(C_+)\subset C_+$ since the Darcy velocity flows outward the injection well. On the contrary, for a production cell $C_-$, we have that $C_- \subset F_{[-\dtDisc,0]}(C_-)$. This indicates that for an efficient application of the GEM scheme, the MMOC component should be implemented on regions near the injection cells once the concentration $c$ is almost constant in these regions. This happens after some time $T_+$ when the injection cells $C_+$ are almost filled up, i.e. $c\approx 1$ in $C_+$.  Before this, we should implement a pure ELLAM scheme. Hence, we start by defining $\alpha=1$ over $\O$ for all $n$ such that $t^{(n)}\le T_+$, the time where the injection cells are filled up; typically, $T_{+}\approx 1 \text{ to } 1.5$ years. Note that $T_+$ can actually be found during the simulation, by checking if the concentration is almost constant in and around the injection cells or not. To be specific, $T_+$ is determined to be the time such that $|\Pi_{\discC} c^{(n)}-\Pi_{\discC}c^{(n+1)}|<\epsilon$ in the cells surrounding the injection well(s) and the well(s) themselves. For our numerical tests, we take $\epsilon=10^{-4}$. For $n$ such that $t^{(n)}>T_+$, and assuming for simplicity one injection cell $C_+$ and one production cell $C_-$, a possible choice is
\be \label{eq:alpha-ini}
\alpha(\x)=\begin{cases}
1 & \quad \mbox{ if $|\x - C_+| \geq |\x - C_-|$ } \\
0 & \quad \mbox{ otherwise. }
\end{cases}
\ee

Here, $|\x - C_+|$ and $|\x -C_-|$ denote the distance between $\x$ and the center of the cells $C_+$ and $C_-$, respectively. This tells us to use ELLAM for regions far from the injection well, and MMOC otherwise. In particular, the choice of $\alpha =1$ near the injection cell and $\alpha = 0$ away from the injection cell ensures that the second condition for mass conservation in Remark \ref{rem.ex:ELLAM-MMOC} is satisfied around the traceback region of $C_+$. Moreover, for regions that are not in the traceback of $C_+$, we have that $\divg\Daru=0$, and hence mass conservation is attained by MMOC c.f. Remarks \ref{rem:comp.ELLAM.MMOC} and \ref{rem:massConsMMOC}. In case of multiple injection and production wells, the same rule can be applied by taking $\alpha(\x)=1$ if the closest well to $\x$ is a production well, and $\alpha(\x)=0$ if the closest well to $\x$ is an injection well.

\subsection{Adaptation of local mass conserving adjustments to the miscible flow model} \label{sec:local-mass-bal-model}
In this section, we write the algorithm proposed in Section \ref{sec:local-mass-bal} in the context of the miscible flow model, for the GDM--ELLAM and GEM scheme. Since the GDM--MMOC does not achieve a global mass balance, it does not make sense to post-process the data and enforce local mass balance for each cell in the mesh. For the miscible flow model \eqref{eq:model}, the velocity field is divergence free in most regions, except for the injection and production wells (see \eqref{eq:press}). Hence, for most cells $K\in\mesh$, we know that $|F_{-\dtDisc}(K)|_\phi = |K|_\phi$.

\subsubsection{GDM--ELLAM}

For the GDM--ELLAM, all cells are tracked backward. Hence, the cells $K$ for which $|F_{-\dtDisc}(K)|$ cannot be determined exactly are: the production cells, and the cells $K \neq C_+$ that track into the injection cells (fully or partially). We thus formulate an approximate local volume constraint for the cells $K \neq C_+$ that track into the injection cells. Due to Liouville's Theorem, the exact volume of the traceback region of an injection cell $C_+$ is deduced to be $|F_{-\dtDisc}(C_+)|_\phi= e^{-\beta}|C_+|_\phi$, where 
\begin{equation}\nonumber
\beta= \dfrac{\int_{C_+} q(t^{(n+1)})}{\int_{C_+} \phi} \dtDisc.
\end{equation}
We compute an approximate trace-forward region $\widetilde{C}_+$ of $C_+$ by forming a polygon with vertices and edge points of $C_+$ tracked forward in time. In practice, compared to the other cells $K$ in the mesh, we track more points along the edges of $C_+$ in order to obtain a good enough approximation of $\widetilde{C}_+$. In particular, if, on average, $n$ points are tracked along the edges of a cell $K \in \mesh$, then we found that $4n+1$ is an appropriate number of points to track along the edges of $C_+$. We then set the approximate local volume constraint for the cells $K\neq C_+$ that track into $C_+$ to be
\begin{equation}\label{eq:approx-loc-vol-cons}
|\widetilde{K}|_\phi = |K|_\phi - |K\cap \widetilde{C}_+|_\phi + \dfrac{|K\cap\widetilde{C}_+|_\phi}{|\widetilde{C}_+ \setminus C_+|_\phi} (1-e^{-\beta}) |C_+|_\phi.
\end{equation} 
Essentially, \eqref{eq:approx-loc-vol-cons} may be interpreted in the following manner: the region $K\cap \widetilde{C}_+$ is the part of $K$ that is expected to track back into the injection cell $C_+$. Hence, since the other part of $K$ stays in the divergence free region, its volume $|K|_\phi - |K\cap \widetilde{C}_+|_\phi$ remains unchanged. The changed volume is then approximated as a ratio of the part in $C_+$ that will be tracked out of itself (i.e. $(1-e^{-\beta}) |C_+|_\phi$). 

We can then adapt the algorithm proposed in Section \ref{sec:local-mass-bal} to the miscible flow model. We note that the local approximation \eqref{eq:approx-loc-vol-cons} is valid only if we have a good approximation of the trace-forward region of the injection cells, and hence we must track more points on the injection cells. As an additional fix, since we expect the injection cells to be eventually filled with the injected fluid, we set $c=1$ on these cells. Assuming that \eqref{eq:approx-loc-vol-cons} gives a good enough approximation, the local volume constraint for the production cells should then be satisfied approximately (since we have global mass conservation).

\subsubsection{GEM}
For the GEM scheme, the only regions with unknown volumes are the trace-forward of the injection cells $F_{\dtDisc}(C_+)$ (for the MMOC component) and the traceback of the production cells $F_{-\dtDisc}(C_-)$ (for the ELLAM component). If the initial approximations made by tracking back the vertices and edges points for the two cells $C_-$ and $C_+$ are good enough, then each of the approximations to the quantities $|F_{-\dtDisc}(C_-)|_\phi$ and $|F_{\dtDisc}(C_+)|_\phi$ should be very close to their actual values. Here, it is acceptable to track more points than what we track on average for the other cells $K\in\mesh$ to ensure that such an accuracy is obtained, since the global impact in terms of cost is minimal. In particular, one clear advantage of the GEM scheme comes from the absence of the approximations \eqref{eq:approx-loc-vol-cons}, which tells us that GEM is expected to give a better local mass conservation property compared to the GDM--ELLAM.

	\begin{remark}[Convergence tests]
		For this miscible displacement flow model, very few analytical solutions exist, and those that are known have a velocity field $\Daru$ such that $\divg(\mathbf{u})$ is  non-zero almost everywhere in the domain (see, e.g., \cite[Section 3.1]{SR21-example}), which means that the wells (corresponding to $q^+$ and $q^-$) are also completely diffused in the domain. Also, when implementing an ELLAM scheme for these test cases, we do not encounter steep back-tracked functions as in Remark \ref{rem:ELLAM.steep}. As a consequence, the most efficient implementation of the  GEM scheme boils down to taking $\alpha=1$ everywhere in \eqref{GSconc-GEM}, corresponding to a pure ELLAM scheme. Due to this, we do not provide convergence tests for the GEM scheme, as it will yield the same results as ELLAM, and would not give any interesting insight on the efficiency of the GEM scheme in preserving global and local mass balances -- its very \emph{raison d'\^etre}. 
\end{remark}               
\subsection{Numerical Results} \label{sec: num-results}

In this section, we compare numerical results obtained from the GEM scheme to those obtained from  GDM--ELLAM. Both of these schemes employ the post-processing technique outlined in Section \ref{sec:local-mass-bal-model} to achieve local mass balance. For completeness, we will also present a comparison with GDM--MMOC, but without the local adjustments. As with the GEM scheme, for the GDM--MMOC, an ELLAM scheme is first implemented for the first few time steps, when $t^{(n)}\leq T_+$, after which, a pure MMOC scheme is implemented, i.e. for $t^{(n)}>T_+$, we take $\alpha=0$ over $\O$ in \eqref{GSconc-GEM}.  Here, we use the HMM scheme \cite{dro-10-uni} for discretising the diffusive terms. This HMM corresponds to a certain choice of the gradient discretisation $\discC$ and $\discP$, see \cite[Chapter 13]{gdm}.
The numerical simulations are performed under the following standard data (see, e.g., \cite{WLELQ-00}):
\begin{enumerate}
	\item $\O=(0,1000) \times (0,1000) \text{ ft}^{2}$, 
	\item injection well at $(1000,1000)$ and production well at $(0,0)$, both with flow rate of $30 \text{ft}^{2}/\text{day}$,
	\item  constant porosity $\phi=0.1$ and constant permeability tensor $\mathbf{K}=80\mathbf{I}$ md,
	\item oil viscosity $\mu(0)=1.0$ cp and mobility ratio $M=41$,
	\item $\phi d_{m}=0.0 \text{ft}^{2}/\text{day}$, $\phi d_{l}=5.0\text{ft}$, and $\phi d_{t}=0.5 \text{ft}$.
\end{enumerate}

\begin{figure}[h]
	\centering
	\begin{tabular}{c@{\hspace*{2em}}c}
		\includegraphics[width=0.4\textwidth]{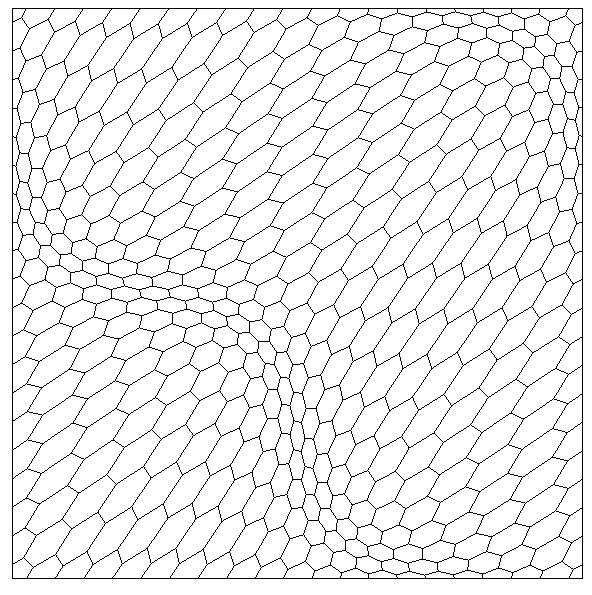} & 		\includegraphics[width=0.4\textwidth]{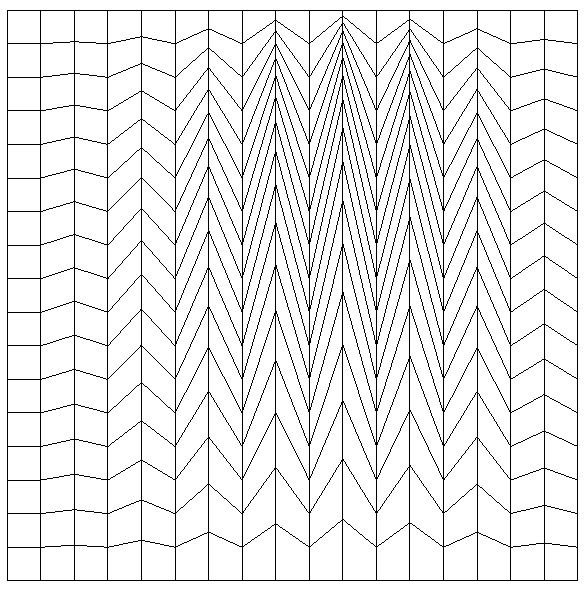}\\
	\end{tabular}
	\caption{ Mesh types(left: hexahedral; right: Kershaw).}
	\label{HKmeshes}
\end{figure}

For the time discretisation, we take a constant time step of $\Delta t =36$ days. The simulations are run on Cartesian  meshes (square cells of dimension $62.5 \times 62.5$ ft), hexahedral meshes, and on Kershaw meshes as described in \cite{HH08} (see Fig. \ref{HKmeshes}). 

\subsubsection{Cartesian meshes}
\begin{figure}[h]
	\begin{tabular}{cc}
		\includegraphics[width=0.4\linewidth]{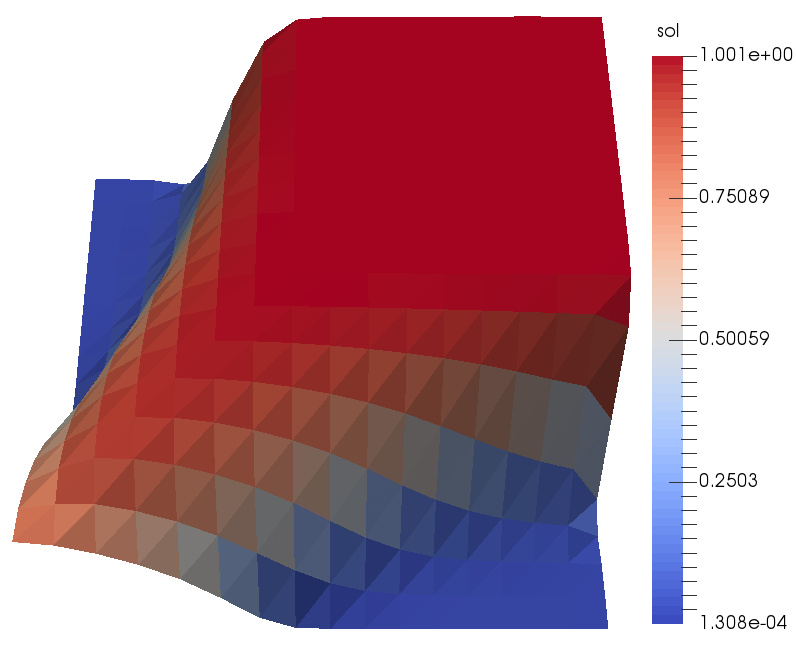} & 		\includegraphics[width=0.4\linewidth]{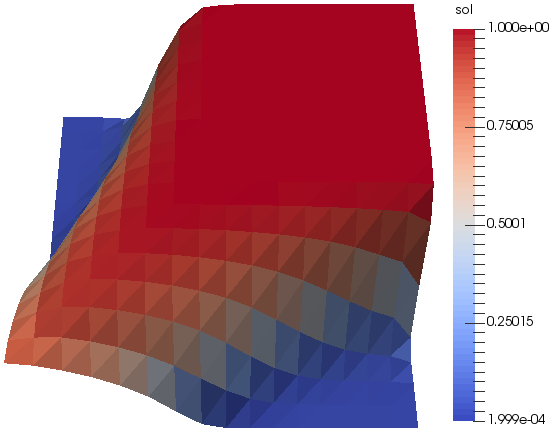}\\
	\end{tabular}
	\caption{concentration profile at $t=10$ years, Cartesian mesh (left: HMM--ELLAM, right: HMM--MMOC).}
	\label{fig:Cart}
\end{figure}

\begin{figure}[h]
	\begin{tabular}{cc}
			\includegraphics[width=0.4\linewidth]{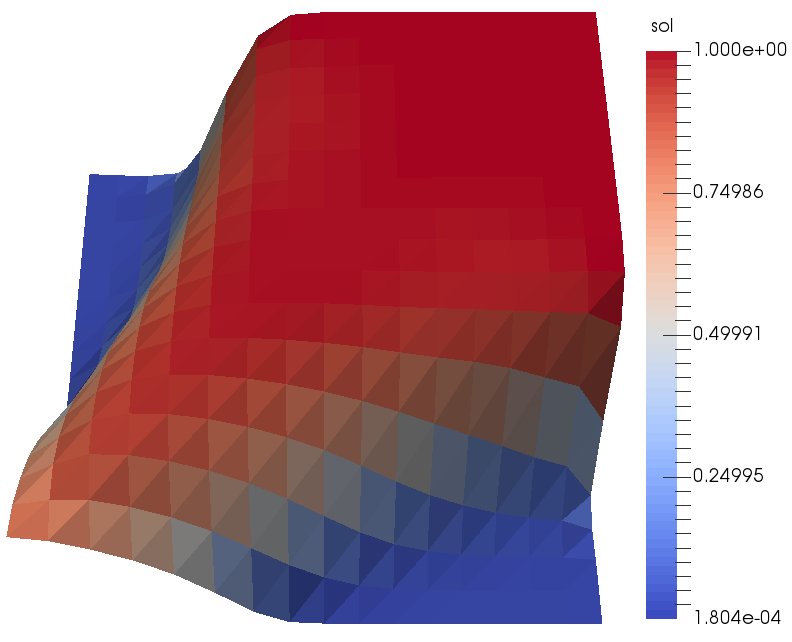} &
		\includegraphics[width=0.45\linewidth]{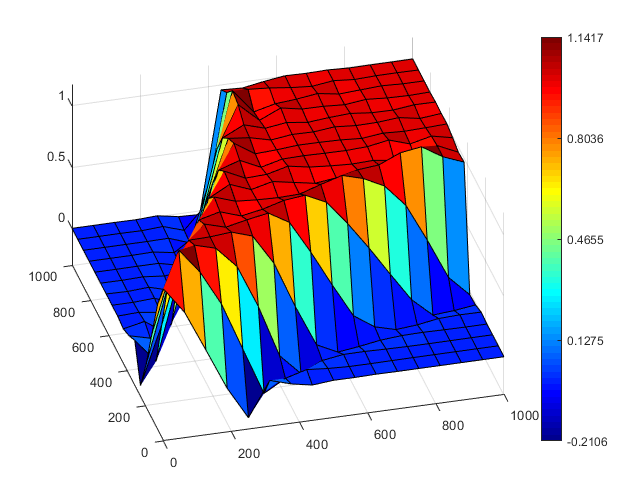}\\
	\end{tabular}

	\caption{concentration profile at $t=10$ years, Cartesian mesh (left: HMM--GEM, right: MFEM--ELLAM).}
	\label{fig:GEM-Cart}
\end{figure}

To compare the numerical solutions, we start by presenting the concentration profiles on a Cartesian mesh at $t=10$ years obtained through HMM--ELLAM and the HMM--MMOC schemes in Figure \ref{fig:Cart}. This is followed by a solution obtained by a HMM--GEM scheme (GEM scheme using the HMM gradient discretisation, with advection components computed as described in \eqref{GSconc-ELLAM-MMOC-pcwise-const}) in Figure \ref{fig:GEM-Cart}. Also, for this test case, we present a numerical result obtained from an MFEM--ELLAM scheme. This numerical output from MFEM--ELLAM was obtained by a straight application of the MFEM--ELLAM algorithm as presented in \cite{WLELQ-00}, with several hundred of quadrature points per cell around the injection well, without any post-processing or tweaks. One of the authors of \cite{WLELQ-00} has informed us that a post-processing was implemented near the injection wells in order to obtain a better concentration profile than those we find in Figure \ref{fig:GEM-Cart}, right. However, due to issues on intellectual property rights, this post-processing technique was not revealed, and hence our tests were not able to account for this. These are accompanied by Table \ref{tab:Cart-GEM}, which presents some important features, such as the number of points tracked along each edge, overshoots/undershoots, $e_{\rm mass}^{(N)}$, and the approximate amount of oil $\frac{1}{|\O|_\phi}\int_\O \phi \Pi_\discC c^{(N)}$ recovered after 10 years. Here, $e_{\rm mass}^{(N)}$ refers to the accumulated mass balance error (percentage) obtained over all the time steps, i.e. 
\[
e_{\rm mass}^{(N)}=\dfrac{\Big|\displaystyle\int_{\O}\phi \Pi_{\discC}c^{(N)}-\int_{\O}\phi \Pi_{\discC}c^{(0)} -\sum_{n=0}^{N-1}\dtDisc\int_{\O} (q^+-\Pi_{\discC} c q^-)^{(n,w)}\Big|}{\Big|\displaystyle\int_{\O}\phi \Pi_{\discC}c^{(0)} +\sum_{n=0}^{N-1}\dtDisc\int_{\O} (q^+-\Pi_{\discC} c q^-)^{(n,w)}\Big|}.
\]
In practice, we have $e_{\rm mass}^{(N)}=\max_{n=1,\dots,N} e_{\rm mass}^{(n)}$, as the mass balance error accumulates over each time step (in the numerical tests, no compensation is observed).

 Upon comparing the concentration profiles, we see that for all of the HMM based schemes, the overshoot is very low, with the maximum overshoot being less than $0.2\%$. However, it can be noted in Table \ref{tab:Cart-GEM} that ELLAM's $0.18\%$ overshoot is much larger, by a factor of almost 20, than those of the MMOC and GEM, even on a very simple square mesh. For the MFEM--ELLAM, the overshoots were at 14.17\%, which is much larger than those observed from the other schemes. Also, we note here the presence of severe undershoots, of around 21\%, near the diagonal. Aside from the overshoots, there are no noticeable differences between the concentration profiles for the HMM--ELLAM and the HMM--GEM. On the other hand, we note that the HMM--MMOC scheme introduces some artificial diffusion along the diagonal, which slightly smears the expected fingering effect.
 \begin{table}[h!]  \centering 		
 	\caption{Comparison between HMM--ELLAM, HMM--MMOC and HMM--GEM schemes, Cartesian mesh}
 	\begin{tabular}{|l|*{4}{c|}}\hline
 		&\makebox{points per edge}&\makebox{overshoot}&\makebox{$e_{\rm mass}^{(N)}$}&\makebox{recovery}\\\hline
 		HMM--ELLAM & 1 & $1.11\%$& 0.19\% & 70.09\% \\\hline
 		HMM--ELLAM & 3 & $0.18\%$& 0.21\% & 69.76\% \\\hline
 		HMM--MMOC & 1 & $<0.01\%$ & 5.60\% & 71.97\% \\\hline
 		HMM--MMOC & 3 & $<0.01\%$ & 2.80\% & 69.94\% \\\hline
 		HMM--GEM & 1 & $<0.01\%$& 2.35\%  &68.44\% \\\hline
 		HMM--GEM & 3 & $<0.01\%$& 0.85\%  &69.14\% \\\hline
 	\end{tabular}
 	\label{tab:Cart-GEM}
 \end{table}

 Next, upon comparing the approximate amount of oil recovered after 10 years, the 68.44\% to 69.14\% obtained for the HMM--GEM scheme is comparable to the amount from the HMM--ELLAM, which ranges from 69.76\% to 70.09\%. The HMM--MMOC, on the other hand, provides an overestimate of the oil recovered when the tracked cells are approximated only by vertices and edge midpoints, due to the excess diffusion it introduces along the diagonal becoming more prominent. When 3 points are tracked along each edge, the amount of oil recovered for HMM--MMOC is almost the same as that for HMM--ELLAM and HMM--GEM. We note that the approximate amount of oil recovered by MFEM--ELLAM is only at 57.39\%, which is possibly caused by the undershoots on the region near the diagonal.

Lastly, we compare the mass balance errors. Since the recovery for MFEM--ELLAM is only at 57.39\%, we expect the mass balance errors to be quite large. Now, we focus on the mass balance errors for the HMM based schemes, obtained once we track 3 or more points along each edge. The error obtained from the GEM (0.85\%) is much better than the one from MMOC (2.80\%), and close to that obtained from ELLAM (0.21\%). These results agree with the analysis provided in Sections \ref{subsec mass-bal-MMOC} and \ref{subsec mass-bal-ELLAM-MMOC}, due to the fact that the MMOC will fail to conserve mass as soon as the fluid starts invading the production well (which translates to $|\Pi_{\discC} c^{(n)}-\Pi_{\discC}c^{(n+1)}|$ being large on $F_{[-\dtDisc,0]}(C_-)$). 

\subsubsection{Hexahedral meshes}

We then compare the numerical results on hexahedral meshes. Unlike the regular square cells for the Cartesian type meshes, the cells for the hexahedral meshes are irregular. To implement the HMM--ELLAM scheme on these types of cells, the proper amount of points to track along each edge is determined by the cell regularity parameter \cite{CD17-HMM-ELLAM-complete}, defined for each cell $K\in\mesh$ to be
\be \nonumber
m_{K\text{reg}} := \dfrac{\text{diam}(K)^2}{|K|}.
\ee
$\lceil \log_{2}(m_{K \text{reg}})\rceil$  points are then tracked along each edge of cell $K$. For hexahedral and Kershaw meshes, we have $\max_{K\in\mesh}(\lceil \log_{2}(m_{K \text{reg}})\rceil)=3$ and $6$, respectively. We note however, that due to the adjustments described in Section \ref{sec:local-mass-bal-model} which need to be implemented in order to achieve local mass conservation, $2 \lceil \log_{2}(m_{K \text{reg}})\rceil + 1 $ points would need to be tracked along the edge of each cell for hexahedral type meshes. In the following figures, three tests are performed for the HMM--ELLAM: the first of which does not involve any adjustment to achieve local mass conservation, followed by an adjustment only on the cells that are not involved with either the injection or production cells (i.e. the algorithm in Section \ref{sec:local-mass-bal-model} without \eqref{eq:approx-loc-vol-cons}), and finally an adjustment based on the full algorithm in Section \ref{sec:local-mass-bal-model}.
\begin{figure}[h]
	\begin{tabular}{cc}
		\includegraphics[width=0.4\linewidth]{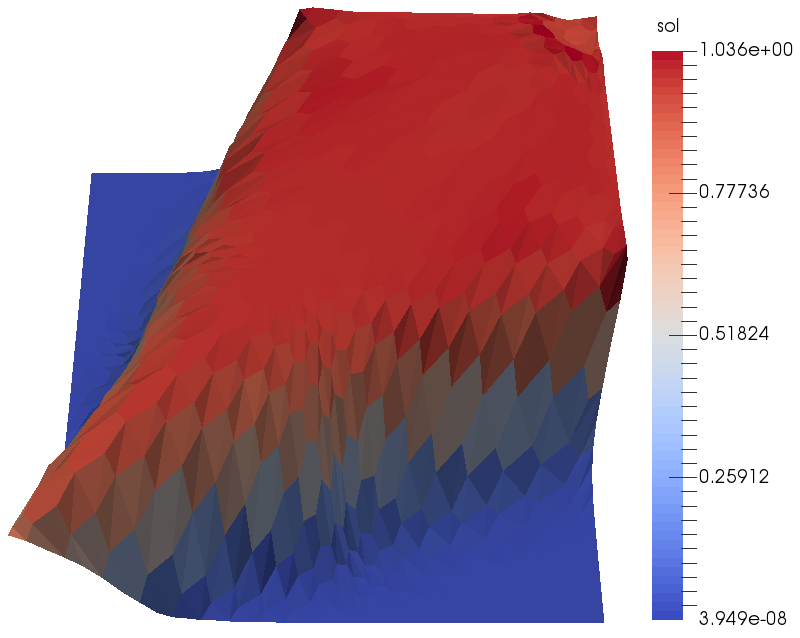} & 		\includegraphics[width=0.4\linewidth]{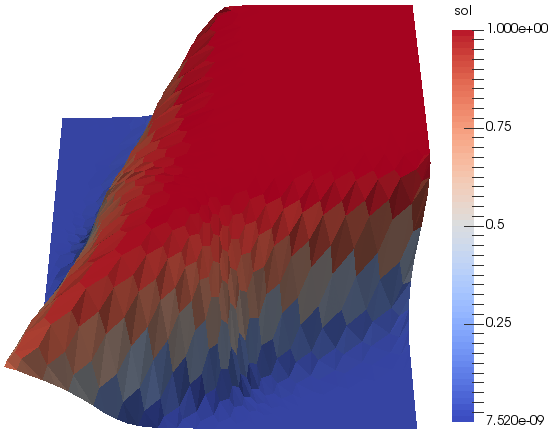}
	\end{tabular}
	\caption{concentration profile at $t=10$ years, hexahedral mesh (left: HMM--ELLAM (no adjustments), right: HMM--MMOC).}
	\label{fig:Hexa}
\end{figure}
\begin{figure}[h]
	\begin{tabular}{cc}
		\includegraphics[width=0.4\linewidth]{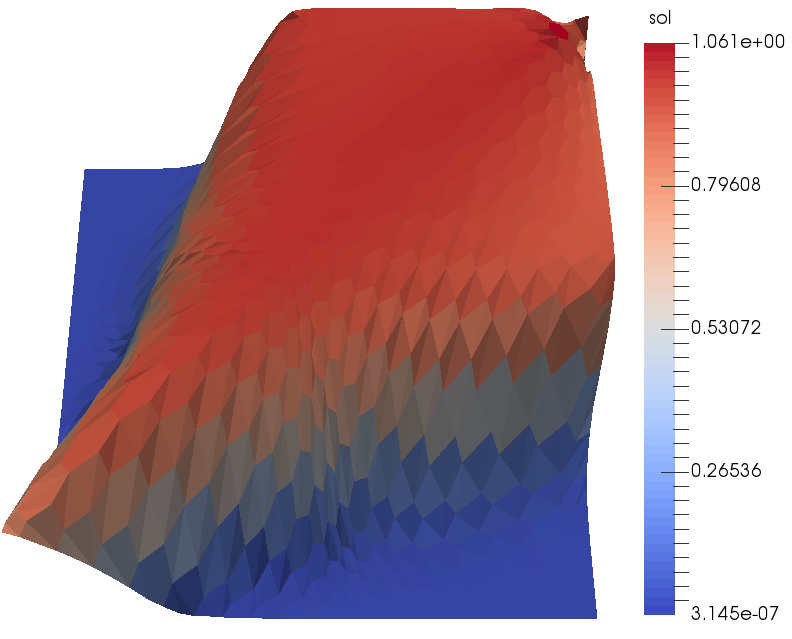} & 		\includegraphics[width=0.4\linewidth]{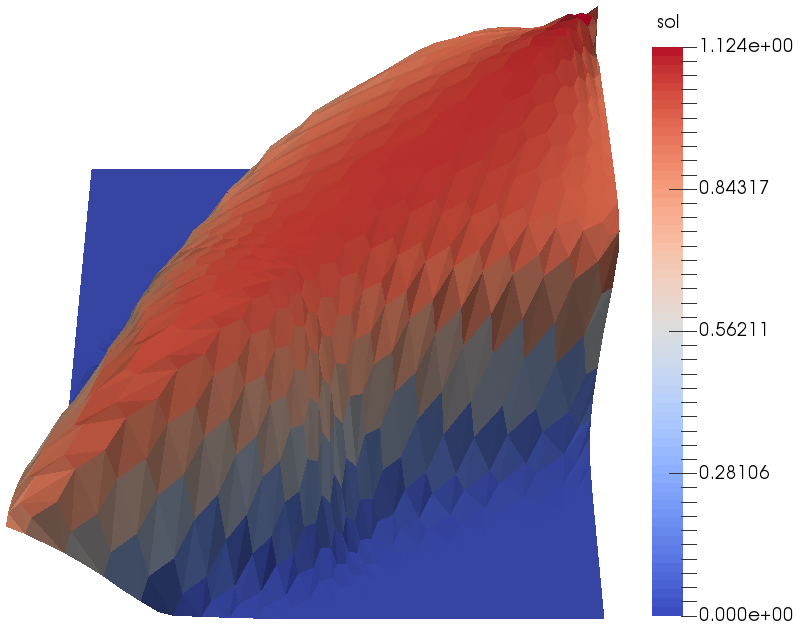}
	\end{tabular}
	\caption{concentration profile at $t=10$ years, hexahedral mesh, HMM--ELLAM with adjustments (left: without \eqref{eq:approx-loc-vol-cons} , right: including \eqref{eq:approx-loc-vol-cons}).}
	\label{fig:Hexa-ELLAM-adj}
\end{figure}
\begin{figure}[h]
\centering
		\includegraphics[width=0.4\linewidth]{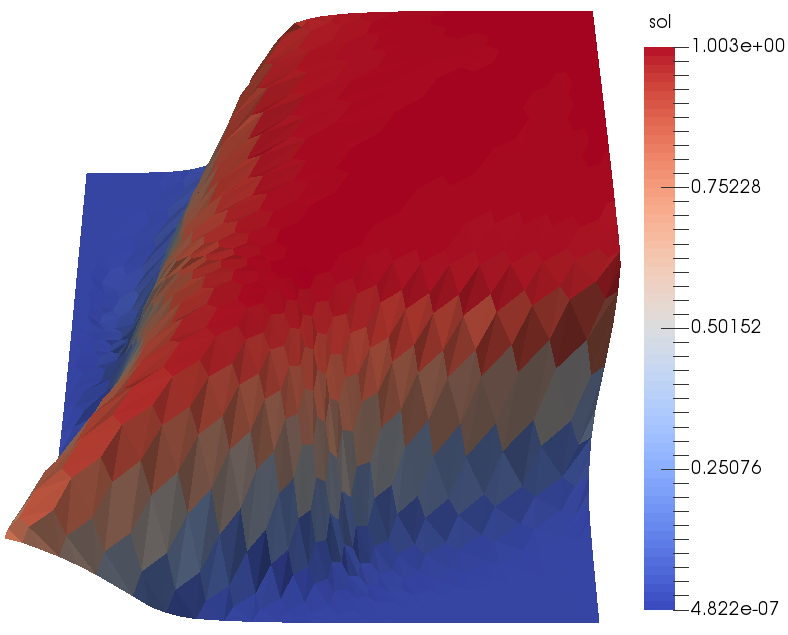}
	\caption{concentration profile at $t=10$ years, hexahedral mesh, GEM scheme.}
	\label{fig:GEM-Hexa}
\end{figure}
\begin{table}[h!]  \centering 		
	\caption{Comparison between HMM--ELLAM, HMM--MMOC and HMM--GEM scheme, hexahedral mesh, $\Delta t = 18$ days}
	\begin{tabular}{|l|*{4}{c|}}\hline
		&\makebox{points per edge}&\makebox{overshoot}&\makebox{$e_{\rm mass}$}&\makebox{recovery}\\\hline
		
		HMM--ELLAM & $\lceil \log_{2}(m_{K \text{reg}})\rceil$ & 3.65\% & 0.62\% & 62.50\% \\
		(no adjustment)  & & & & \\\hline 
		HMM--ELLAM  & $2\lceil \log_{2}(m_{K \text{reg}})\rceil+1$ & 6.14\% & 0.31\% & 62.81\% \\
		(adjustment without \eqref{eq:approx-loc-vol-cons}) & & & & \\\hline 
		HMM--ELLAM  & $2\lceil \log_{2}(m_{K \text{reg}})\rceil+1$ & 12.42\% & 0.74\% & 57.04\% \\
		(adjustment with \eqref{eq:approx-loc-vol-cons}) & & & & \\\hline 
		HMM--MMOC & $\lceil \log_{2}(m_{K \text{reg}})\rceil$ & $<0.01\%$ & 1.82\% & 61.43\% \\\hline
		HMM--GEM & $2\lceil \log_{2}(m_{K \text{reg}})\rceil+1$ & $0.34\%$& 0.74\% & 64.13\% \\\hline
	\end{tabular}
	\label{tab:Hexa-GEM}
\end{table}

Upon looking at the results in Figures \ref{fig:Hexa}-\ref{fig:GEM-Hexa} and Table \ref{tab:Hexa-GEM}, it is noticeable that the solution from HMM--ELLAM has a large discrepancy and overshoot near the injection well. This is due to the fact that the algorithm in Section \ref{sec:local-mass-bal-model} fails to converge for hexahedral meshes. We note here that such a failure of the adjustment algorithm was not noticed on Cartesian meshes, whether in the tests conducted above or (with a different approach to the adjustment) in \cite{AW11-stability-monotonicity-implementation}. Our tests on hexahedral meshes demonstrate here the difficulty of designing a robust algorithm to locally adjust the mass balance for ELLAM. As a matter of fact, it has been pointed out in \cite{D16-Opti-meshCorr} that there is no guarantee that adjustments for ELLAM type schemes, such as those in \cite{AW11-stability-monotonicity-implementation}, will terminate or yield a valid mesh configuration. The same issue could happen to our proposed adjustment applied on the GEM method but, as seen in the tests above, this adjustment seems to be more robust.

There are two possible explanations for why the algorithm in Section \ref{sec:local-mass-bal-model} fails to converge for hexahedral meshes: Firstly, compared to the Cartesian mesh, the volume of the injection cell, and each of the cells around it (around 700 to 1000 square units), is much smaller than the volume of the other cells in the mesh (on average, 2000 square units). Since these cells are already small to begin with, tracking them backwards will lead to traceback regions which are much smaller, and hence will be more prone to errors. Taking $\Delta t = 36$ days, the expected volume of the traceback region of the injection cell $e^{-\alpha}|C_+| \approx 1e-04$ square units, which is only around 1e-08\% of the total volume. This gives us an idea that taking a smaller time step might be able to mitigate these errors; however, even taking a much smaller time step of $\Delta t = 1$ day, there is still no convergence.  Due to this, we conclude that it is caused by the second possibility, i.e. for hexahedral meshes, \eqref{eq:approx-loc-vol-cons} does not give a good approximation of the local volume constraint for the cells tracked back into the injection cell. This can be seen more clearly upon comparing Figure \ref{fig:Hexa-ELLAM-adj} (left and right). Without enforcing \eqref{eq:approx-loc-vol-cons}, the solution seems to behave better. Actually, upon implementation of \eqref{eq:approx-loc-vol-cons}, the local volume constraint is satisfied on the cells which are not tracked back either into injection or into production cells. Having eliminated the local mass balance errors from these cells, they accumulate onto the cells near the injection cell. Since \eqref{eq:approx-loc-vol-cons} does not give a good approximation, the accumulated error around this region does not spread and cancel out properly, and hence severely distorts the quality of the numerical solution. Contrary to the HMM--ELLAM, due to the absence of \eqref{eq:approx-loc-vol-cons}, this problem is not as severe with the GEM scheme, and can be solved by taking a smaller time step of $\Delta t = 18$. 
  
With the exception of the case for which ELLAM is adjusted with the local volume constraint \eqref{eq:approx-loc-vol-cons}, the amount of oil recovered from all three schemes are comparable, as they are within 2\% of each other. Upon comparing the mass balance errors, we note that the HMM--GEM (0.74\%) outperforms the HMM--MMOC (1.82\%), and is in the same range as the best HMM--ELLAM implementation (0.31\%). This example shows that, on non-Cartesian meshes, due to the absence of \eqref{eq:approx-loc-vol-cons}, the HMM--GEM is able to provide a better-looking solution, with reduced overshoots and acceptable mass conservation properties, compared to the HMM--ELLAM method. Moreover, the HMM--GEM conserves mass better than the HMM--MMOC. 

\subsubsection{Kershaw meshes}\label{subsubsec: num-Kershaw}
Finally, we compare the numerical results on a much more challenging mesh, the Kershaw mesh.
\begin{figure}[h]
	\begin{tabular}{cc}
		\includegraphics[width=0.4\linewidth]{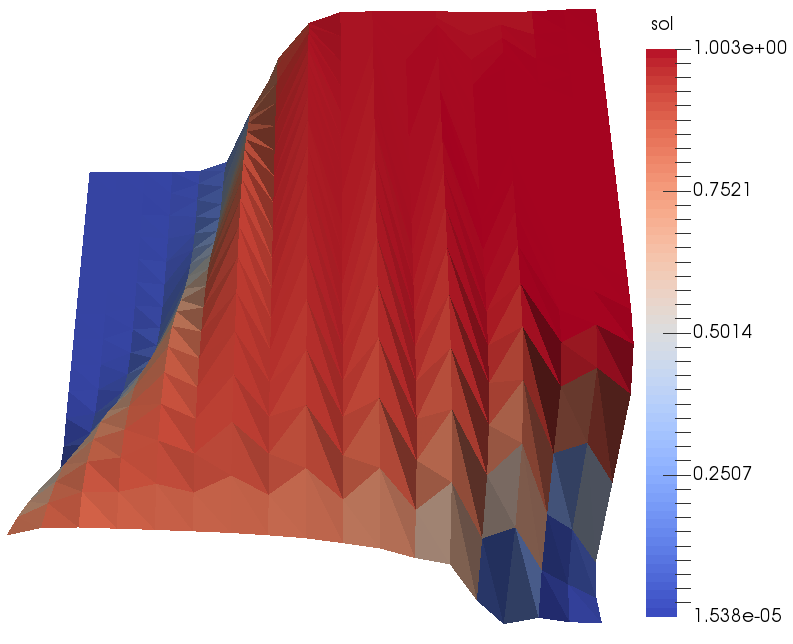} & 		\includegraphics[width=0.4\linewidth]{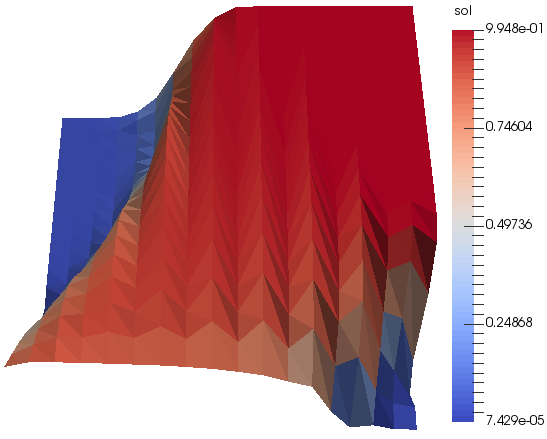}
	\end{tabular}
	\caption{concentration profile at $t=10$ years, Kershaw mesh (left: HMM--ELLAM, right: HMM--MMOC).}
	\label{fig:Ker}
\end{figure}

 \begin{figure}[h]
	\centering
	\includegraphics[width=0.4\linewidth]{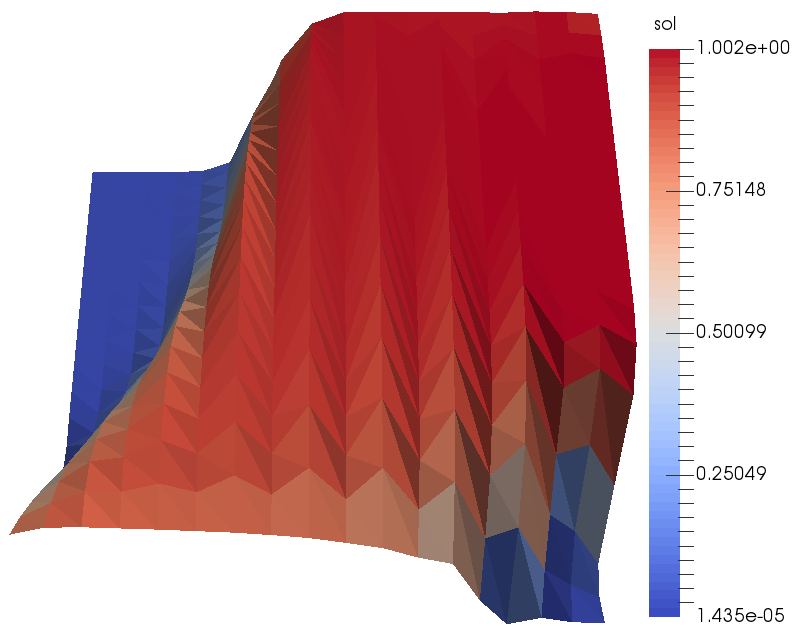}
	\caption{concentration profile at $t=10$ years, Kershaw mesh, HMM--GEM.}
	\label{fig:GEM-Ker-ini}
\end{figure}

\begin{table}[h!]  \centering 		
	\caption{Comparison between HMM--ELLAM, HMM--MMOC and HMM--GEM scheme, Kershaw mesh}
	\begin{tabular}{|l|*{4}{c|}}\hline
		&\makebox{points per edge}&\makebox{overshoot}&\makebox{$e_{\rm mass}$}&\makebox{recovery}\\\hline
		HMM--ELLAM & $\lceil \log_{2}(m_{K \text{reg}})\rceil$ & $0.28\%$  & 0.38\% & 72.63\% \\\hline 
		HMM--MMOC & $\lceil \log_{2}(m_{K \text{reg}})\rceil$ & $0\%$  & $4.28\%$ & 73.21\% \\\hline
		HMM--GEM  & $\lceil \log_{2}(m_{K \text{reg}})\rceil$ & $0.20\%$ & 0.23\% & 72.52\% \\\hline
	\end{tabular}
	\label{tab:Ker-GEM1}
\end{table}

 As with the Cartesian mesh test case, no significant difference can be observed between the numerical solutions obtained from HMM--ELLAM and HMM--GEM. Also, as expected, the mass balance error for the HMM--MMOC is quite large. Notably, the numerical solution on Kershaw type meshes is skewed towards the lower right corner. This is expected due to the fact that the numerical fluxes for HMM schemes on this type of mesh are prone to grid effects, as explained in \cite{CD17-HMM-ELLAM-complete}. 

\medskip
To summarise the previous tests, the HMM--GEM exhibits a slightly better volume conservation property than the HMM--ELLAM. This is due to the local volume constraint for the HMM--ELLAM being inexact in the sense that it depended on \eqref{eq:approx-loc-vol-cons} to give a good approximation; whereas for the HMM--GEM scheme, exact local volume constraints were imposed. In particular, on non-Cartesian meshes with small or mildly distorted cells, HMM--GEM can control local volume constraints more easily and more accurately than HMM--ELLAM, as was seen in the test case on hexahedral meshes. Hence, with the choice of $\alpha$ driven by the discussion in Section \ref{subsec mass-bal-ELLAM-MMOC}, GEM achieves both a good preservation of the physical bounds on $c$, and of mass conservation.

\section{Convergence result for GDM--MMOC scheme}\label{sec:conv.anal}

In this section, we present a convergence result for the GDM--MMOC scheme. To simplify the exposition, this convergence is stated for the scalar advection--reaction--diffusion model \eqref{eq:advection}; at the expense of heavier notations, a convergence result of GDM--MMOC could be stated for the coupled model \eqref{eq:model}, as in \cite{CDL17-convergence-ELLAM} for the GDM--ELLAM method. This convergence result is established using only weak regularity assumptions on the solution (which are satisfied in practical applications),
and not relying on $L^\infty$ bounds (which are impossible to ensure at the discrete level given
the anisotropic diffusion tensors and the general grids used in applications). Hence, the convergence is established under compactness arguments, and error estimates are not available. As a matter of fact, the convergence analysis of numerical approximations of \eqref{eq:model} under weak regularity assumptions has recently received an increasing interest; see, e.g., \cite{CD-07,ckm15} for finite volume
methods and \cite{rivwalk11,GLR17-conv-DG} for discontinuous Galerkin methods. We also note that this convergence analysis assumes a perfect computation of the tracked regions. An analysis which accounts for the approximation in tracked regions, and adjustment strategies, has been performed in \cite{AW10-convergence}, under the assumption that $\Daru \in C^1(Q_T)^d$. Future work will aim to perform a convergence analysis, without this strong regularity assumption.

We assume here \eqref{assump.global}, as well as, denoting by $\mathcal S_d(\R)$ the set of $d\times d$
symmetric matrices,
\begin{equation}\label{hyp:diffTens}
\begin{aligned}
&	\Lambda:Q_T\to \mathcal S_d(\R)\mbox{ is measurable, bounded and uniformly coercive:}\\
& \quad\exists \nu_\Lambda>0\,,\;\exists \alpha_\Lambda>0\mbox{ s.t., for a.e.\ $\x\in\R^d$, a.e. $t\in (0,T)$ and all $\xi\in\R^d$,}\\
& \quad|\Lambda(\x,t)|\leq \nu_\Lambda\quad\text{and}\quad 
	|\Lambda(\x,t)\xi\cdot \xi|\geq \alpha_\Lambda|\xi|^2.
\end{aligned}
\end{equation}
Adding the diffusion term to the MMOC scheme \eqref{conc-MMOC} for the advection--reaction model, we obtain the definition of the GDM--MMOC scheme for \eqref{eq:advection}.
\begin{definition}[GDM--MMOC scheme]
	The GDM--MMOC scheme for \eqref{eq:advection} reads as: find 
	 $(c^{(n)})_{n=0,\ldots,N}\in X_{\discC}^{N+1}$ such that $c^{(0)}=\ICinterp_{\discC} c_{\rm ini}$
	 and, for all $n=0,\ldots,N-1$,
		$c^{(n+1)}$ satisfies
		\begin{equation}\label{conc-MMOC2}
		\begin{aligned}
		\int_{\O} \phi {}&\PiDc c^{(n+1)} \PiDc zd\x -\int_{\O} \phi(\x) (\PiDc c^{(n)})(F_{-\dtDisc}(\x)) \PiDc z(\x) d\x\\
		&+\dtDisc \int_{\O} \Lambda^{(n+1)} \gradDc c^{(n+1)} \cdot \gradDc z d\x\\
		={}&\dtDisc \int_{\O} \bigl[\bigl(f^{(n,w)}-(\Pi_{\discC} c)^{(n,w)} \divg \Daru^{(n+1)}\bigr)\cdot \mathbf{e}\bigr]\PiDc z d\x,\qquad \forall z\in X_\discC,
		\end{aligned}
		\end{equation}
		where $\Lambda^{(n+1)}(\x):=\frac{1}{\dtDisc}\int_{t^{(n)}}^{t^{(n+1)}}\Lambda(\x,s)ds$
\end{definition}

If $c=(c^{(n)})_{n=0,\ldots,N}$ is a solution to the GDM--MMOC scheme, we define $\Pi_\discC c\in L^\infty(Q_T)$ and $\nabla_\discC c\in L^\infty(Q_T)^d$ by
\[
\begin{aligned}
\forall n=0,\ldots,N-1\,,\;\forall t\in (t^{(n)},t^{(n+1)})\,,\;\mbox{ for a.e.\ $\x\in\O$},\\
\Pi_\discC c(\x,t)=\Pi_\discC c^{(n+1)}(\x)\mbox{ and }\nabla_\discC c(\x,t)=\nabla_\discC c^{(n+1)}(\x)
\end{aligned}
\]

The convergence of the GDM--MMOC can be established under most of the assumptions made for the GDM--ELLAM (i.e. \cite[Assumptions \combineln{A}{1},\combineln{A}{3}-\combineln{A}{4}]{CDL17-convergence-ELLAM}). Assumption \ref{hyp:approx.un} makes precise the sense in which the velocities used to construct the flows over each time step (see \eqref{charac}) approximate the given velocity $\Daru$. 
\begin{enumerate}[label=\combineln{B}{\arabic*},leftmargin=2.8em]
%

	\item \label{hyp:approx.un} For each $m\in\N$, approximate velocities $(\Daru_m^{(n+1)})_{n=0,\ldots,N-1}$ are chosen in $L^2(\O)^d$, with $\div\Daru_m^{(n+1)}\in L^\infty(\O)$,
and there are $M_{\div},M_{\rm vel}\ge 0$ such that, for all $m\in\N$ and $n=0,\ldots,N_m-1$,
\[
\Vert \Daru_m^{(n+1)}\Vert_{L^2(\Omega)}\le M_{\rm vel}\quad\mbox{ and }\quad
\Vert\div \Daru_m^{(n+1)}\Vert_{L^\infty(\Omega)}\le M_{\rm div}.
\]
Letting $\bar\Daru_m\in L^\infty(0,T;L^2(\O)^d)$ be defined by $\bar\Daru_m(\x,t)=\Daru_m^{(n+1)}(\x)$ for a.e.\ $\x\in\O$, all $t\in (t_m^{(n)},t_m^{(n+1)})$ and all $n=0,\ldots,N_m-1$, we have, as $m\to\infty$, $\bar\Daru_m\to\Daru$ weakly in $L^2(Q_T)^d$.

%
	
\end{enumerate}

We now define the notion of weak solution to the advection--reaction--diffusion model, and state the convergence result for the GDM--MMOC. The proof of this convergence, that relies on discrete compactness techniques, is similar to the proof of \cite[Theorem 3.3]{CDL17-convergence-ELLAM} and is therefore omitted.

\begin{definition}[Weak solution to the advection--reaction--diffusion equation]
	A function $c\in L^2(0,T;H^1(\O))$ is a weak solution of \eqref{eq:advection} if 
	it satisfies
	\begin{equation}\label{adv.weak}
	\begin{aligned}
	&-\int_\O \phi(\x)c_{\rm ini}(\x)\varphi(\x,0) d\x
	-\int_0^T\int_{\Omega} \phi(\x) c(\x,t)\dfrac{\partial \varphi}{\partial t}(\x,t) d\x dt\\
	&-\int_0^T\int_{\Omega} c(\x,t)\darcyU(\x,t) \cdot \nabla \varphi(\x,t)d\x dt
	+\int_0^T\int_{\Omega} \Lambda\nabla c(\x,t)\cdot \nabla \varphi(\x,t)d\x dt
 \\&=\int_0^T\int_{\Omega} f(c,\x,t)\varphi(\x,t)d\x dt, \quad \forall \varphi \in C^\infty_c(\overline{\O}\times [0,T)).
	\end{aligned}
	\end{equation}
\end{definition}
\begin{theorem}[Convergence of the GDM--MMOC]\label{th:MMOC convergence}
	Under Assumptions \eqref{assump.global}, \eqref{hyp:diffTens}, \ref{hyp:approx.un}, and \cite[Assumptions \combineln{A}{1}, \combineln{A}{3}-\combineln{A}{4}]{CDL17-convergence-ELLAM}, for any $m\in\N$ there is a unique
	$c_m\in X_{\discC_m}^{N_m+1}$ solution of the MMOC scheme \eqref{conc-MMOC} with $
	\discC^T=\discC_m^T$ and $(\Daru^{(n+1)})_{n=0,\ldots,N-1}=(\Daru_m^{(n+1)})_{n=0,\ldots,N_m-1}$.
	Moreover, up to a subsequence as $m\to\infty$,
	\begin{itemize}
		\item $\Pi_{\discC_m} c_m \conv c$ weakly-$*$ in $L^\infty(0,T;L^2(\O))$ and strongly
		in $L^r(0,T;L^2(\O))$ for all $r<\infty$,
		\item $\nabla_{\discC_m}c_m\to \nabla c$ weakly in $L^2(Q_T)^d$,
	\end{itemize}
	where $c$ is a weak solution of \eqref{eq:advection}.
\end{theorem}


\medskip

\thanks{\textbf{Acknowledgement}: This research was supported by the Australian Government through the Australian Research Council's Discovery Projects funding scheme (pro\-ject number DP170100605).
}

\bibliographystyle{abbrv}
\bibliography{char-methods}
\end{document}